\documentclass[11pt,reqno]{amsart}

\usepackage{graphicx,booktabs}
\usepackage{amsmath,amssymb,mathrsfs,cite}
\usepackage{subfigure}
\usepackage{color}

\newtheorem{theorem}{Theorem}[section]

\theoremstyle{definition}

\theoremstyle{remark}
\newtheorem{remark}[theorem]{Remark}
\numberwithin{equation}{section}

\begin{document}

\title[A fast algorithm for the electromagnetic cavity scattering problem]{A
fast algorithm for the electromagnetic scattering from a large
rectangular cavity in three dimensions}

\author{Yanli Chen}
\address{Department of Mathematics, Northeastern University, Shenyang 110819,
China}
\email{chenyanli@mail.neu.edu.cn}

\author{Xue Jiang}
\address{School of Mathematics, Faculty of Science, Beijing University of Technology, Beijing,
100124, China.}
\email{jxue@lsec.cc.ac.cn}

\author{Jun Lai}
\address{School of Mathematical Sciences, Zhejiang University Hangzhou, Zhejiang
310027, China}
\email{laijun6@zju.edu.cn}

\author{Peijun Li}
\address{Department of Mathematics, Purdue University, West Lafayette, Indiana
47907, USA}
\email{lipeijun@math.purdue.edu}




\begin{abstract}
The paper is concerned with the three-dimensional electromagnetic scattering
from a large open rectangular cavity that is embedded in a perfectly
electrically conducting infinite ground plane. By introducing a transparent
boundary condition, the scattering problem is formulated into a boundary value
problem in the bounded cavity. Based on the Fourier expansions of the
electric field, the Maxwell equation is reduced to one-dimensional ordinary
differential equations for the Fourier coefficients. A fast algorithm,
employing the fast Fourier transform and the Gaussian elimination, is
developed to solve the resulting linear system for the cavity which is filled
with either a homogeneous or a layered medium. In addition, a novel scheme is
designed to evaluate rapidly and accurately the Fourier transform of singular
integrals. Numerical experiments are presented for large
cavities to demonstrate the superior performance of the proposed method.
\end{abstract}

\maketitle

\section{Introduction}

The electromagnetic scattering from large cavities has received much attention
in both engineering and mathematical communities due to its significant
industrial and military applications\cite{BGLZ12,BS05,J91,LAG14,DSZ13}. For
instance, the radar cross section (RCS) measures the detectability of a
target by a radar system. In practice, the cavity RCS caused by objects such
as jet engine inlet ducts, exhaust nozzles and cavity-backed antennas can
dominate the total RCS. Therefore, mathematical and computational methods to
accurately predict the cavity RCS are important for the enhancement or reduction
of the total RCS~\cite{BJL1,BJL2}.  Another example is the non-destructive
testing to determine the shape of a cavity embedded in a known object. In these
applications, it has played a crucial role to have an efficient forward solver
for the optimal design problems of reducing or enhancing the cavity RCS and the
inverse problems of determining an unknown cavity.

A variety of numerical methods, including finite difference methods, finite
element methods, the moment methods, boundary element methods, and hybrid
methods, have been developed to solve the open cavity problems
\cite{JV91,VW03,ZMH09,WG03,WW99,HCH07,LMS13,WDS08}, In particular, Bao and Sun
\cite{BS05} proposed a finite difference based fast algorithm for the
two-dimensional electromagnetic scattering from large cavities. In the
algorithm, an FFT-sine transform in the horizontal direction and the Gaussian
elimination along the vertical direction were used to reduce the global system
to a much smaller system imposed only on the open aperture of the cavity. As an
extension of this method, a tensor product finite element method was proposed
in \cite{DSZ13} by employing piecewise polynomials of degree $k\geq 1$ to
approximate the solution space of the cavity problem. In \cite{ZQT11}, a fourth
order finite difference scheme was developed to discrete the cavity scattering
problem in the rectangular domain and to reach a global fourth order
convergence in the whole computational domain by a special treatment on the
boundary condition. Since the resulting linear system obtained from the cavity
problem is usually  indefinite and ill-conditioned, convergence of iterative
methods such as GMRES is very slow. Different kinds of
preconditioners were proposed to accelerate the convergence \cite{BS05,DSZ13,
ZQT11,D11,ZZ19}. On the other hand, a fast direct solver based on hierarchical
matrix factorization technique was used to solve the two-dimensional
electromagnetic scattering from an arbitrarily shaped cavity\cite{LAG14}.  It
was shown that the linear system resulted from the integral equation method can
be solved in nearly linear time. The method was extended to the scattering of
three-dimensional axis-symmetric cavities in \cite{LAI20171}. We refer to
\cite{J02} for the motivation, modeling, computation, as well as related
references on the open cavity scattering problems.

It is worth mentioning that the computation is extremely challenging when the
cavities are large compared to the wavelength of the incident wave because of
the highly oscillatory nature of the fields.  For such a high frequency
scattering problem, it is shown that the ratio of the error by the usual
Galerkin type method and the error of the best approximation tends to infinity
as the wave number increases \cite{BS97, AKS88}. Due to these difficulties,
the discretization by conventional numerical methods becomes very
expensive for the large cavity scattering problems especially in three
dimensions. In this paper, we intend to develop a fast algorithm for solving the
three-dimensional electromagnetic scattering from large rectangular cavities
embedded in an infinite perfectly electrically conducting ground plane.

More specifically, we consider the three-dimensional Maxwell equations along
with the Silver--M\"{u}ller radiation condition imposed at infinity. By using
the dyadic Green's function in the half space, we first derive an exact
transparent boundary condition (TBC) on the open aperture of the
cavity. As a result, the original scattering problem is formulated equivalently
to a boundary value problem of Maxwell equations in a bounded domain. Secondly,
we introduce the Fourier series expansion of the electric field inside the
cavity. By such an expansion, the governing Maxwell equations can be reduced to
one-dimensional ordinary differential equations with respect to the vertical
direction. A second-order finite difference scheme is adopted to solve the
ordinary differential systems. A fast algorithm, based on the fast Fourier
transform in the horizontal directions and the Gaussian elimination along the
vertical direction, is developed to solve the linear system arising from
scattering of large cavities which may be filled with a homogeneous medium or a
vertically layered medium. Moreover, we reduce the global system to a linear
system on the open aperture of the cavity only and design a novel scheme to
evaluate rapidly and accurately the singular integrals appeared in the
transparent boundary condition. Numerical results show that our algorithm is
very efficient in terms of computational cost.

The paper is organized as follows. In Section 2, we describe the problem
formulation of the electromagnetic scattering by a rectangular cavity which is
filled with a homogeneous medium. The governing Maxwell equations along with the
Silver--M\"{u}ller radiation condition are introduced. The TBC is presented to
reduce the unbounded scattering problem to a boundary value problem formulated
in the bounded cavity. The details of the fast algorithm are given in Section
3. Section 4 is devoted to an extension of the fast algorithm to the scattering
of a cavity which is filled with a layered medium. Section 5 proposes an FFT
based efficient algorithm to evaluate the singular integrals arising from the
nonlocal TBC on the open aperture of the cavity. Analysis on the computational
complexity for the fast algorithm is discussed in Section 6. Numerical examples
are presented in Section 7 to demonstrate the performance of the proposed
algorithm. The paper is concluded with some general remarks in Section 8.

\section{Problem formulation}

\begin{figure}
 \includegraphics[width=0.5\textwidth]{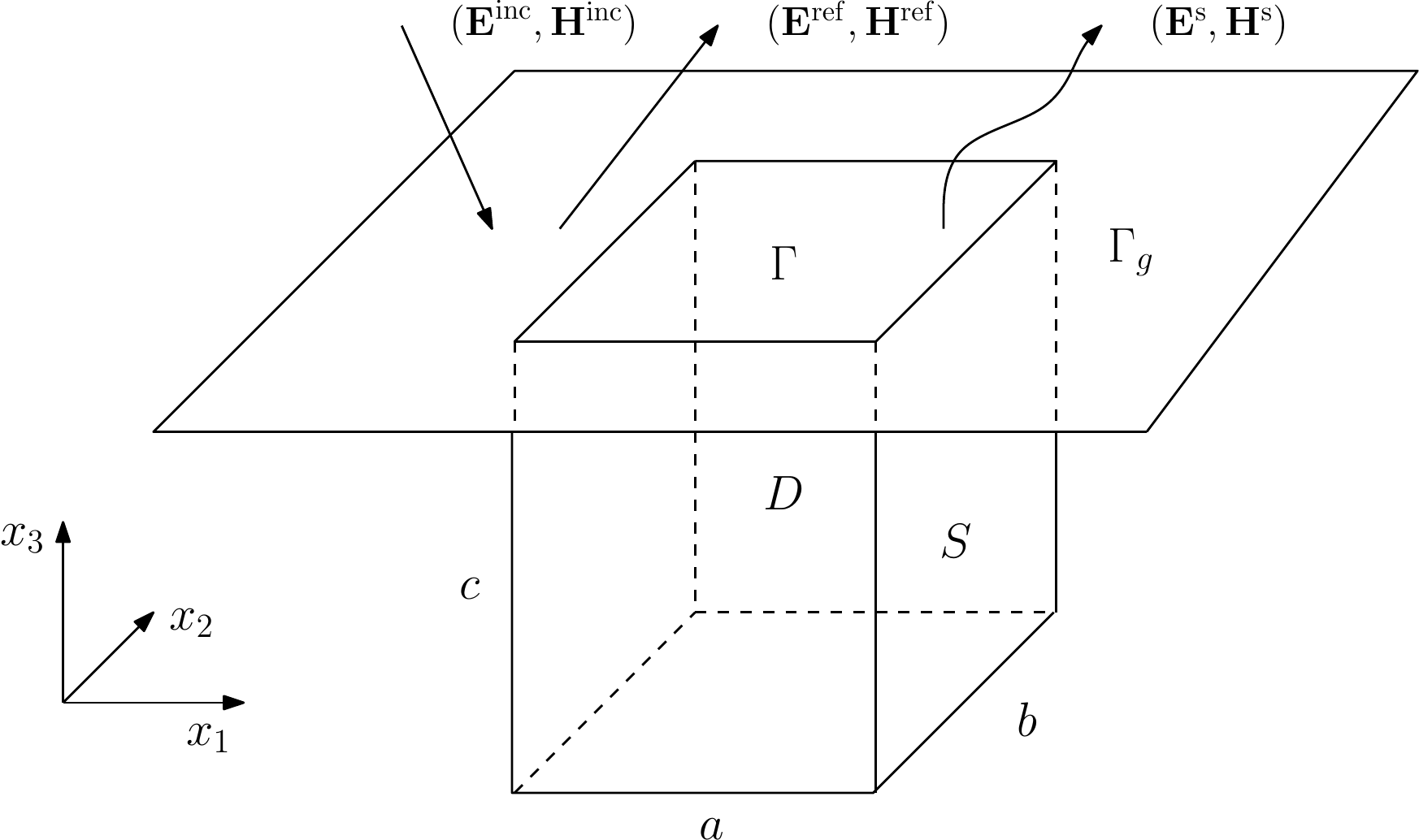}
 \caption{The problem geometry of the electromagnetic scattering by a
rectangular cavity.}
 \label{fig:1}
 \end{figure}

Consider the incidence of a time-harmonic electromagnetic wave on a rectangular
cavity $D\subset \mathbb{R}^3$, which is embedded in the infinite ground plane
$\Gamma_g$. The problem geometry is shown in Figure \ref{fig:1}. The cavity wall
$S$ and the ground plane $\Gamma_g$ are assumed to be perfect electric
conductors. We also assume that the open aperture $\Gamma=[0,a]\times[0,b]$ is
aligned with the ground plane $\Gamma_g$ and the depth of the cavity is $c$. The
half space above the ground and the cavity are assumed to be filled with some
homogeneous material with a constant electric permittivity $\varepsilon_0$ and a
constant magnetic permeability $\mu_0$. Let $B_R^+$ be a half-ball above the
ground plane with hemisphere  $\Gamma_R^+$ as part of the boundary, where the
radius $R$ is large enough so that $\Gamma_R^+$ covers the open aperture
$\Gamma$. It is clear to note that the full boundary of $\partial B_R^+$
consists of the hemisphere $\Gamma_R^+$, the open aperture $\Gamma$, and a part
of the ground plane $\Gamma_g$. Without confusion, we simply denote  $\partial
B_R^+=\Gamma_R^+\cup\Gamma\cup\Gamma_g$.

The total electric and magnetic fields $(\boldsymbol{E},\boldsymbol{H})$
consist of the incident waves $(\boldsymbol{E}^{\rm inc},\boldsymbol{H}^{\rm
inc})$, the reflected waves $(\boldsymbol{E}^{\rm ref},\boldsymbol{H}^{\rm
ref})$ due to the infinite ground plane, and the scattered wave
$(\boldsymbol{E}^{\rm s},\boldsymbol{H}^{\rm s})$ because of the open cavity.
The total fields $\boldsymbol{E}$ and $\boldsymbol{H}$ satisfy Maxwell's
equations $\mathrm{in}~\mathbb{R}_+^3\cup D$:
\begin{equation}\label{Meq}
\nabla\times \boldsymbol{E}=i\omega\mu_0\boldsymbol{H},\quad\nabla \times \boldsymbol{H}=-i\omega \varepsilon_0 \boldsymbol{E},
\end{equation}
where $\omega>0$ is the angular frequency. Since the ground plane and the cavity
wall are perfect conductors, we have
\begin{equation}\label{PEC}
\nu\times\boldsymbol{E}=0 \quad {\rm on}~\Gamma_g\cup S,
\end{equation}
where $\nu$ is the unit normal vector on $\Gamma_g$ and $S$.

The incident electromagnetic plane waves $(\boldsymbol{E}^{\rm
inc},\boldsymbol{H}^{\rm inc})$ are given as
\begin{equation*}
  \boldsymbol{E}^{\rm inc}=\boldsymbol{p}e^{i \boldsymbol{q}\cdot\boldsymbol{x}},\quad\boldsymbol{H}^{\rm inc}=\boldsymbol{s}e^{i\boldsymbol{q}\cdot\boldsymbol{x}},\quad \boldsymbol{s}=\frac{\boldsymbol{q}\times\boldsymbol{p}}{\omega\mu_0},\quad\boldsymbol{p}\cdot\boldsymbol{q}=0,
\end{equation*}
where $\boldsymbol x=(x_1, x_2, x_3)\in\mathbb
R^3$, $\boldsymbol{p}=(p_1,p_2,p_3)$ and $\boldsymbol{s}=(s_1,s_2,s_3)$ are the
polarization vectors, $\boldsymbol{q}=(\alpha_1,\alpha_2,-\beta)$ with $\beta\ge
0$ is the propagation direction vector.
It is easy to verify that the incident electromagnetic fields
$(\boldsymbol{E}^{\rm inc},\boldsymbol{H}^{\rm inc})$ satisfy the Maxwell equation
\eqref{Meq} in $\mathbb{R}_+^3$.

Due to the infinite ground plane, the reflected fields $(\boldsymbol{E}^{\rm ref},\boldsymbol{H}^{\rm ref})$ can be explicitly written as
\begin{equation*}
  \boldsymbol{E}^{\rm ref}=\boldsymbol{p}^*e^{i \boldsymbol{q}^*\cdot\boldsymbol{x}},\quad\boldsymbol{H}^{\rm ref}=\boldsymbol{s}^*e^{i\boldsymbol{q}^*\cdot\boldsymbol{x}},\quad \boldsymbol{s}^*=\frac{\boldsymbol{q}^*\times\boldsymbol{p}^*}{\omega\mu_0},\quad\boldsymbol{p}^*\cdot\boldsymbol{q}^*=0,
\end{equation*}
where $\boldsymbol{p}^*=(-p_1,-p_2,p_3)$ and $\boldsymbol{q}^*=(\alpha_1,\alpha_2,\beta)$. Evidently, the reflected fields $(\boldsymbol{E}^{\rm ref},\boldsymbol{H}^{\rm ref})$ also satisfy the Maxwell equation  \eqref{Meq} in $\mathbb{R}_+^3$. In particular, the following homogeneous Dirichlet boundary condition is satisfied for the incident and reflected electric fields on the  ground plane:
\begin{equation*}
  \nu\times(\boldsymbol{E}^{\rm inc}+\boldsymbol{E}^{\rm ref})=0\quad {\rm on}~ \Gamma_g.
\end{equation*}

It follows from \eqref{Meq} and the incident and reflected electromagnetic
fields that the scattered electromagnetic fields $(\boldsymbol E^{\rm s},
\boldsymbol H^{\rm s})$ also satisfy the Maxwell equation
\begin{equation}\label{EHs}
\nabla\times \boldsymbol{E}^{\rm s}=i\omega\mu_0
\boldsymbol{H}^{\rm s},\quad\nabla \times \boldsymbol{H}^{\rm s}=-i\omega
\varepsilon_0 \boldsymbol{E}^{\rm s},\quad \boldsymbol{x}\in\mathbb{R}_+^3,
\end{equation}
and the homogeneous Dirichlet boundary condition
\begin{equation}\label{Esb}
  \nu\times\boldsymbol{E}^s=0\quad {\rm on}~ \Gamma_g.
\end{equation}
In addition, the scattered field $(\boldsymbol{E}^{\rm s},\boldsymbol{H}^{\rm
s})$ are required to satisfy the Silver--M\"{u}ller radiation
condition:
\begin{equation}\label{SMRC}
\sqrt{\varepsilon_0}\boldsymbol{E}^{\rm s}-\sqrt{\mu_0}\boldsymbol{H}^{\rm
s}\times\hat{\boldsymbol{x}}=o(|\boldsymbol{x}|^{-1}),\quad
|\boldsymbol{x}|\rightarrow\infty,
\end{equation}
where $\hat{\boldsymbol{x}}=\boldsymbol{x}/|\boldsymbol{x}|$. By eliminating
the scattered magnetic field in \eqref{EHs}, the scattered electric field
satisfies
\begin{equation}\label{Esweq}
\nabla\times(\nabla\times \boldsymbol{E}^{\rm s})-\kappa_0^2\boldsymbol{E}^{\rm s}=0\quad \mathrm{in}~\mathbb{R}_+^3,
\end{equation}
where $\kappa_0=\omega\sqrt{\varepsilon_0\mu_0}$ is the wavenumber.

In order to derive a transparent boundary condition on the open aperture
$\Gamma$, we introduce the half-space dyadic Green's
function $\bar{\bar{G}}_e(\boldsymbol{x},\boldsymbol{y})$, which is given by
\begin{equation}\label{HDGF}
\bar{\bar{G}}_e(\boldsymbol{x},\boldsymbol{y})=\bar{\bar{G}}_0(\boldsymbol{x},\boldsymbol{y})-\bar{\bar{G}}_0(\boldsymbol{x},\boldsymbol{y}_i)+2\hat{\boldsymbol{z}}\hat{\boldsymbol{z}}g(\boldsymbol{x},\boldsymbol{y}_i),
\end{equation}
where
\begin{equation}\label{FDGF}
\bar{\bar{G}}_0(\boldsymbol{x},\boldsymbol{y})=\Big{(}\bar{\bar{I}}-\frac{1}{\kappa_0^2}\nabla_{\boldsymbol{x}}\nabla_{\boldsymbol{y}}\Big{)}g(\boldsymbol{x},\boldsymbol{y}),
\end{equation}
is the free space dyadic Green's function,
$\bar{\bar{I}}=\hat{\boldsymbol{x}}\hat{\boldsymbol{x}}+\hat{\boldsymbol{y}}\hat
{\boldsymbol{y}}+\hat{\boldsymbol{z}}\hat{\boldsymbol{z}}$ is the
$3\times 3$ identity matrix, and
\begin{equation}\label{GF}
  g(\boldsymbol{x},\boldsymbol{y})=\frac{e^{i\kappa_0|\boldsymbol{x}-\boldsymbol{y}|}}{4\pi|\boldsymbol{x}-\boldsymbol{y}|},
\end{equation}
is the free space Green's function for the three-dimensional Helmholtz
equation. Here
$\boldsymbol{y}_i=y_1\hat{\boldsymbol{x}}+y_2\hat{\boldsymbol{y}}-y_3\hat{
\boldsymbol{z}}$ denotes the image point of
$\boldsymbol{y}=y_1\hat{\boldsymbol{x}}+y_2\hat{\boldsymbol{y}}+y_3\hat{
\boldsymbol{z}}$, and $\hat{\boldsymbol x}, \hat{\boldsymbol y},
\hat{\boldsymbol z}$ are the unit vectors in the $x_1, x_2, x_3$ axis,
respectively.

The half-space dyadic Green's function satisfies the Maxwell
equation
\begin{equation}\label{Geq}
\nabla\times\big{(}\nabla\times \bar{\bar{G}}_e(\boldsymbol{x},\boldsymbol{y})\big{)}-\kappa_0^2\bar{\bar{G}}_e(\boldsymbol{x},\boldsymbol{y})=\bar{\bar{I}}\delta(\boldsymbol{x}-\boldsymbol{y})\quad {\rm in}~\mathbb{R}_+^3,
\end{equation}
and the Dirichlet boundary condition
\begin{equation}\label{Gb}
\nu\times\bar{\bar{G}}_e(\boldsymbol{x},\boldsymbol{y})=0\quad {\rm on}~
\Gamma_g\cup \Gamma,
\end{equation}
where $\delta$ is the Dirac delta function. Furthermore, the half-space dyadic
Green's function satisfies the Silver--M\"{u}ller radiation condition.

Next, we present the transparent boundary condition. Multiplying both sides of
\eqref{Esweq} by the half-space dyadic Green's function and integrating over
$B_R^+$, we obtain
\begin{equation*}
  \int_{B_R^+}\Big{(}\big{(}\nabla_{\boldsymbol{x}}\times\nabla_{\boldsymbol{x}}\times\boldsymbol{E}^{\rm s}(\boldsymbol{x})\big{)}\cdot\bar{\bar{G}}_e(\boldsymbol{x},\boldsymbol{y})-\kappa_0^2\boldsymbol{E}^{\rm s}(\boldsymbol{x})\cdot\bar{\bar{G}}_e(\boldsymbol{x},\boldsymbol{y})\Big{)}{\rm d}\boldsymbol{x}=0.
\end{equation*}
It follows from the second vector Green's theorem that
\begin{align}
\nonumber   &\int_{B_R^+}\boldsymbol{E}^{\rm
s}(\boldsymbol{x})\cdot\Big{(}\nabla_{\boldsymbol{x}}\times\nabla_{\boldsymbol{x
}}\times\bar{\bar{G}}_e(\boldsymbol{x},\boldsymbol{y})-\kappa_0^2\bar{\bar{G}}
_e(\boldsymbol{x},\boldsymbol{y})\Big{)}{\rm d}\boldsymbol{x}  \\\label{EG}
   &=-\int_{\Gamma_R^+\cup\Gamma\cup
\Gamma_g}\Big{(}\big{(}\nu\times\boldsymbol{E}^{\rm
s}(\boldsymbol{x})\big{)}\cdot\big{(}\nabla_{\boldsymbol{x}}\times\bar{\bar{G}}
_e(\boldsymbol{x},\boldsymbol{y})\big{)}-\big{(}\nu\times\bar{\bar{G}}
_e(\boldsymbol{x},\boldsymbol{y})\big{)}\cdot\big{(}\nabla_{\boldsymbol{x}}
\times\boldsymbol{E}^{\rm s}(\boldsymbol{x})\big{)}\Big{)}{\rm
d}s_{\boldsymbol{x}}.
\end{align}
Since the scattered field $\boldsymbol{E}^{\rm s}(\boldsymbol{x})$ and the
half-space dyadic Green's function satisfy the Silver--M\"{u}ller
radiation condition, we get
\begin{equation}\label{GR+}
\int_{\Gamma_R^+}\Big{(}\big{(}\nu\times\boldsymbol{E}^{\rm s}(\boldsymbol{x})\big{)}\cdot\big{(}\nabla_{\boldsymbol{x}}\times\bar{\bar{G}}_e(\boldsymbol{x},\boldsymbol{y})\big{)}-\big{(}\nu\times\bar{\bar{G}}_e(\boldsymbol{x},\boldsymbol{y})\big{)}\cdot\big{(}\nabla_{\boldsymbol{x}}\times\boldsymbol{E}^{\rm s}(\boldsymbol{x})\big{)}\Big{)}{\rm d}s_{\boldsymbol{x}}=0.
\end{equation}
Combining \eqref{Esb} and \eqref{Gb} gives
\begin{equation}\label{Gg}
\int_{\Gamma_g}\Big{(}\big{(}\nu\times\boldsymbol{E}^{\rm
s}(\boldsymbol{x})\big{)}\cdot\big{(}\nabla_{\boldsymbol{x}}\times\bar{\bar{G}}
_e(\boldsymbol{x},\boldsymbol{y})\big{)}-\big{(}\nu\times\bar{\bar{G}}
_e(\boldsymbol{x},\boldsymbol{y})\big{)}\cdot\big{(}\nabla_{\boldsymbol{x}}
\times\boldsymbol{E}^{\rm s}(\boldsymbol{x})\big{)}\Big{)}{\rm
d}s_{\boldsymbol{x}}=0
\end{equation}
and
\begin{equation}\label{Ga}
\int_{\Gamma}\big{(}\nu\times\bar{\bar{G}}_e(\boldsymbol{x},\boldsymbol{y})\big{)}\cdot\big{(}\nabla_{\boldsymbol{x}}\times\boldsymbol{E}^{\rm s}(\boldsymbol{x})\big{)}{\rm d}s_{\boldsymbol{x}}=0.
\end{equation}
Using \eqref{EG}--\eqref{Ga} yields
\begin{align}
\nonumber   &\int_{B_R^+}\boldsymbol{E}^{\rm
s}(\boldsymbol{x})\cdot\Big{(}\nabla_{\boldsymbol{x}}\times\nabla_{\boldsymbol{x
}}\times\cdot\bar{\bar{G}}_e(\boldsymbol{x},\boldsymbol{y})-\kappa_0^2\bar{\bar{
G}}_e(\boldsymbol{x},\boldsymbol{y})\Big{)}{\rm d}\boldsymbol{x}
\\\label{EsGa1}
   &= -\int_{\Gamma}\big{(}\nu\times\boldsymbol{E}^{\rm
s}(\boldsymbol{x})\big{)}\cdot\big{(}\nabla_{\boldsymbol{x}}\times\bar{\bar{G}}
_e(\boldsymbol{x},\boldsymbol{y})\big{)}{\rm d}s_{\boldsymbol{x}}.
\end{align}
Substituting \eqref{Geq} into \eqref{EsGa1} and switching variables
$\boldsymbol{x}$ and $\boldsymbol{y}$, we get
\begin{equation*}
\boldsymbol{E}^{\rm s}(\boldsymbol{x})=-\int_{\Gamma}\big{(}\nu\times\boldsymbol{E}^{\rm s}(\boldsymbol{y})\big{)}\cdot\big{(}\nabla_{\boldsymbol{y}}\times\bar{\bar{G}}_e(\boldsymbol{x},\boldsymbol{y})\big{)}{\rm d}s_{\boldsymbol{y}}.
\end{equation*}
Noting $\nu=-\hat{\boldsymbol{z}}$ gives
\begin{equation*}
\boldsymbol{E}^{\rm s}(\boldsymbol{x})=\int_{\Gamma}\big{(}\hat{\boldsymbol{z}}\times\boldsymbol{E}^{\rm s}(\boldsymbol{y})\big{)}\cdot\big{(}\nabla_{\boldsymbol{y}}\times\bar{\bar{G}}_e(\boldsymbol{x},\boldsymbol{y})\big{)}{\rm d}s_{\boldsymbol{y}}.
\end{equation*}
It follows from $\boldsymbol{E}^{\rm s}=\boldsymbol{E}-\boldsymbol{E}^{\rm inc}-\boldsymbol{E}^{\rm ref}$ and $\hat{\boldsymbol{z}}\times(\boldsymbol{E}^{\rm inc}+\boldsymbol{E}^{\rm ref})=0$ on $\Gamma$ that
\begin{equation}\label{EGa1}
\boldsymbol{E}=\boldsymbol{E}^{\rm inc}+\boldsymbol{E}^{\rm ref}+\int_{\Gamma}\big{(}\hat{\boldsymbol{z}}\times\boldsymbol{E}(\boldsymbol{y})\big{)}\cdot\big{(}\nabla_{\boldsymbol{y}}\times\bar{\bar{G}}_e(\boldsymbol{x},\boldsymbol{y})\big{)}{\rm d}s_{\boldsymbol{y}}.
\end{equation}
Substituting \eqref{HDGF} into \eqref{EGa1}, we obtain
\begin{equation}\label{EGa2}
\boldsymbol{E}=\boldsymbol{E}^{\rm inc}+\boldsymbol{E}^{\rm ref}+2\int_{\Gamma}\big{(}\hat{\boldsymbol{z}}\times\boldsymbol{E}(\boldsymbol{y})\big{)}\cdot\big{(}\nabla_{\boldsymbol{y}}\times\bar{\bar{G}}_0(\boldsymbol{x},\boldsymbol{y})\big{)}{\rm d}s_{\boldsymbol{y}}.
\end{equation}
Taking curl on the both sides of \eqref{EGa2} yields
\begin{equation}\label{ENU}
\nabla_{\boldsymbol{x}}\times\boldsymbol{E}=\nabla_{\boldsymbol{x}}\times\boldsymbol{E}^{\rm inc}+\nabla_{\boldsymbol{x}}\times\boldsymbol{E}^{\rm ref}-2\kappa_0^2\int_{\Gamma}\big{(}\hat{\boldsymbol{z}}\times\boldsymbol{E}(\boldsymbol{y})\big{)}\cdot\bar{\bar{G}}_0(\boldsymbol{x},\boldsymbol{y}){\rm d}s_{\boldsymbol{y}}.
\end{equation}
Substituting \eqref{FDGF} into \eqref{ENU}, we get
\begin{align*}
\big{(}\nabla_{\boldsymbol{x}}\times\boldsymbol{E}\big{)}&=\nabla_{\boldsymbol{
x}}\times\boldsymbol{E}^{\rm inc}+\nabla_{\boldsymbol{x}}\times\boldsymbol{E}^{\rm ref}
-2\kappa_0^2\int\limits_{\Gamma}\big{(}\hat{\boldsymbol{z}}\times\boldsymbol{E}
(\boldsymbol{y})\big{)}
g(\boldsymbol{x},\boldsymbol{y}){\rm d}s_{\boldsymbol{y}}\\
&\quad +2\Big{(}\nabla_{\boldsymbol{x}}\int_{\Gamma}\big{(}\hat{\boldsymbol{z}}
\times\boldsymbol{E}(\boldsymbol{y})\big{)}\cdot\big{(}\nabla_{\boldsymbol{y}}
g(\boldsymbol{x},\boldsymbol{y})\big{)}{\rm d}s_{\boldsymbol{y}}\Big{)}.
\end{align*}
For a continuous differential function $u$ defined in a neighborhood of $\Gamma$, define the surface gradient on $\Gamma$ by
\begin{equation*}
  \nabla_{\Gamma}u=(\nu\times\nabla u)\times \nu.
\end{equation*}
Moreover, we have the decomposition
\begin{equation}\label{decom}
\nabla u= \nabla_{\Gamma}u+\frac{\partial u}{\partial\nu} \nu,
\end{equation}
where $\frac{\partial u}{\partial\nu}$ is the normal derivative on $\Gamma$. Let
$\boldsymbol{v}$ be a tangent vector on $\Gamma$, then we have
\begin{equation}\label{onGa}
\int_{\Gamma}u{\rm div}_{\Gamma}\boldsymbol{v}{\rm
d}s=-\int_{\Gamma}\nabla_{\Gamma}u\cdot\boldsymbol{v}{\rm d}s.
\end{equation}
Using \eqref{decom}--\eqref{onGa} and taking the limit $x_3\rightarrow 0+$, we
obtain the following transparent boundary condition (TBC):
\begin{equation}\label{ENU1}
\hat{\boldsymbol{z}}\times\big{(}\nabla_{\boldsymbol{x}}\times\boldsymbol{E}\big
{)}=\mathscr{T}(\boldsymbol{E}) +\boldsymbol g\quad {\rm on}~ \Gamma,
\end{equation}
where $\boldsymbol g=\hat{\boldsymbol{z}}\times\big{(}\nabla_{
\boldsymbol{x}} \times\boldsymbol{E}^{\rm
inc}\big{)}+\hat{\boldsymbol{z}}\times\big{(}\nabla_{
\boldsymbol{x}}\times\boldsymbol{E}^{\rm ref}\big{)} $ and
\begin{align*}
\mathscr{T}(\boldsymbol{E})=
-2\kappa_0^2\hat{\boldsymbol{z}}\times\int_{\Gamma}\big{(}\hat{\boldsymbol{z}}
\times\boldsymbol{E}(\boldsymbol{y})\big{)} g(\boldsymbol{x},\boldsymbol{y}){\rm
d}s_{\boldsymbol{y}}-2\hat{\boldsymbol{z}}\times\Big{(}\nabla_{\boldsymbol{x}}
\int_{\Gamma}{ \rm
div}_{\Gamma}\big{(}\hat{\boldsymbol{z}}\times\boldsymbol{E}(\boldsymbol{y}
)\big {)}g(\boldsymbol{x},\boldsymbol{y}){\rm d}s_{\boldsymbol{y}}\Big{)}.
\end{align*}
Then, by eliminating the magnetic field in \eqref{Meq} and using the TBC
\eqref{ENU1}, the scattering problem \eqref{Meq}--\eqref{PEC} can be reduced to
an equivalent boundary value problem in the cavity $D$:
\begin{equation}\label{beq1}
\left\{\begin{aligned}
&\nabla\times(\nabla\times \boldsymbol{E})-\kappa_0^2\boldsymbol{E}=0 \quad&& \mathrm{in}~D,\\
&\nu\times\boldsymbol{E}=0 \quad&& {\rm on}~S,\\
&\hat{\boldsymbol{z}}\times\big{(}\nabla_{\boldsymbol{x}}\times\boldsymbol{E}
\big{)}=\mathscr{T}(\boldsymbol{E})+\boldsymbol g\quad&& {\rm on}~ \Gamma.
\end{aligned}\right.
\end{equation}

\section{Discretization and fast algorithm}

In this section, we present the numerical discretization to the Maxwell
equation and the TBC, and a fast algorithm for the resulting system.

Let $\boldsymbol{E}=(E_1,E_2,E_3)$. On the plane surfaces $x_1=0$ and $x_1=a$,
the unit outward normal vectors are $(-1,0,0)$ and $(1,0,0)$, respectively.
Using the boundary condition in \eqref{beq1}, we get the homogeneous Dirichlet
boundary condition for $E_2$ and $E_3$:
\begin{equation}\label{E23D}
E_2(0,x_2,x_3)=E_2(a,x_2,x_3),\quad
E_3(0,x_2,x_3)=E_3(a,x_2,x_3).
\end{equation}
Recall the divergence free condition on the surface:
\begin{equation*}
  \nabla\cdot \boldsymbol{E}=\partial_{x_1}E_1+\partial_{x_2}E_2+\partial_{x_3}E_3=0,
\end{equation*}
which, together with \eqref{E23D}, implies the homogeneous Neumann boundary condition for $E_1$:
\begin{equation}\label{E1N}
\partial_{x_1}E_1(0,x_2,x_3)=\partial_{x_1}E_1(a,x_2,x_3)=0.
\end{equation}
Similarly, on the plane surfaces $x_2=0$ and $x_2=b$, the unit outward normal
vectors are $(0,-1,0)$ and $(0,1,0)$, respectively. Using the boundary condition
in \eqref{beq1}, we have the homogeneous Dirichlet boundary condition for $E_1$ and
$E_3$:
\begin{equation}\label{E13D}
E_1(x_1,0,x_3)=E_1(x_1,b,x_3),\quad
E_3(x_1,0,,x_3)=E_3(x_1,b,x_3).
\end{equation}
Using \eqref{E13D} and the divergence free condition again gives the homogeneous Neumann boundary condition for $E_2$:
\begin{equation}\label{E2N}
\partial_{x_2}E_2(x_1,0,x_3)=\partial_{x_2}E_2(x_1,b,x_3)=0.
\end{equation}
By the boundary conditions \eqref{E23D}--\eqref{E2N}, it is easy to show that
$E_j, j=1,2,3$ admits the following Fourier series expansions:
\begin{equation}\label{E123F}
  \left\{\begin{aligned}
  E_1(x_1,x_2,x_3)&=\sum\limits_{k\in
\mathbb{N}^2}E_1^{(k)}(x_3)\cos\Big{(}\frac{k_1\pi
x_1}{a}\Big{)}\sin\Big{(}\frac{k_2\pi x_2}{b}\Big{)},\\
  E_2(x_1,x_2,x_3)&=\sum\limits_{k\in
\mathbb{N}^2}E_2^{(k)}(x_3)\sin\Big{(}\frac{k_1\pi x_1
}{a}\Big{)}\cos\Big{(}\frac{k_2\pi x_2}{b}\Big{)},\\
  E_3(x_1,x_2,x_3)&=\sum\limits_{k\in
\mathbb{N}^2}E_3^{(k)}(x_3)\sin\Big{(}\frac{k_1\pi x_1
}{a}\Big{)}\sin\Big{(}\frac{k_2\pi x_2}{b}\Big{)},
  \end{aligned}\right.
  \end{equation}
where $k=(k_1,k_2)\in \mathbb{N}^2$.

By the vector identity $\nabla\times(\nabla\times \boldsymbol{E})=-\Delta
\boldsymbol{E}+\nabla(\nabla\cdot \boldsymbol{E})$ and the divergence free
condition $\nabla\cdot \boldsymbol{E}=0$, the Maxwell equation in \eqref{beq1}
can be reduced to the vector Helmholtz equation
 \begin{equation}\label{Heq}
  \Delta \boldsymbol{E}+\kappa_0^2\boldsymbol{E}=0 \quad \mathrm{in}~ D.
\end{equation}
Using the boundary condition in \eqref{beq1} and the divergence free condition
on the plane surfaces $x_3=-c$, we get the homogeneous Dirichlet boundary
condition
\begin{equation}\label{E12D}
  E_1(x_1,x_2,-c)=E_2(x_1,x_2,-c)=0
\end{equation}
and  the homogeneous Neumann boundary condition
\begin{equation}\label{E3N}
\partial_{x_3}E_3(x_1,x_2,-c)=0.
\end{equation}

Substituting \eqref{E123F} into \eqref{Heq}--\eqref{E3N}, we may get the second
order ordinary differential equations for the Fourier coefficients $E_l^{(m,n)},
l=1, 2$:
\begin{equation}\label{ode1x3}
  \left\{\begin{aligned}
  &\frac{{\rm d}^2} {{\rm d}x_3^2} E_l^{(m,n)}(x_3)+\Big{(}\kappa_0^2-\big{(}\frac{m\pi}{a}\big{)}^2-\big{(}\frac{n\pi}{b}\big{)}^2\Big{)}E_l^{(m,n)}(x_3)=0,\quad x_3\in(-c,0),\\
  &E_l^{(m,n)}(-c)=0,
  \end{aligned}\right.
  \end{equation}
 where $(m,n)\in\mathbb{N}^2_l$,  and the second
order ordinary differential equations for the Fourier coefficients
$E_3^{(m,n)}$:
\begin{equation}\label{ode3x3}
  \left\{\begin{aligned}
  &\frac{{\rm d}^2}{{\rm d}x_3^2}E_3^{(m,n)}(x_3)+\Big{(}\kappa_0^2-\big{(}\frac{m\pi}{a}\big{)}^2-\big{(}\frac{n\pi}{b}\big{)}^2\Big{)}E_3^{(m,n)}(x_3)=0,\quad x_3\in(-c,0),\\
  &\frac{\rm d}{{\rm d}x_3}E_3^{(m,n)}(-c)=0.
  \end{aligned}\right.
  \end{equation}
where $(m,n)\in\mathbb{N}^2_3$. Here
$\mathbb{N}^2_1=\{0,1,2,\cdots,M\}\times\{1,2,\cdots,N\}$,
$\mathbb{N}^2_2=\{1,2,\cdots,M\}\times\{0,1,2,\cdots,N\}$ and
$\mathbb{N}^2_3=\{1,2,\cdots,M\}\times\{1,2,\cdots,N\}$, $M$ and $N$ are the
finite truncation numbers of the Fourier series.

Let $\{x_3^j\}_{j=0}^{j=J+1}$ be a set of uniformly distributed grid points of
$[-c,0]$ with $x_3^{j+1}-x_3^j=h$. Let $E_{l,j}^{(m,n)}$ be the finite
difference solution of $E_l^{(m,n)}(x_3), l=1,2,3$ at the point $x_3=x_3^j$. The
discrete finite difference systems for \eqref{ode1x3}--\eqref{ode3x3} are
\begin{equation*}
\left\{\begin{aligned}
&\frac{E_{l,j-1}^{(m,n)}-2E_{l,j}^{(m,n)}+E_{l,j+1}^{(m,n)}}{h^2}+\Big{(}\kappa_0^2-\big{(}\frac{m\pi}{a}\big{)}^2-\big{(}\frac{n\pi}{b}\big{)}^2\Big{)}E_{l,j}^{(m,n)}=0,~ j=1,2,\cdots,J,\\
&E_{l,0}^{(m,n)}=0,
\end{aligned}\right.
\end{equation*}
and
\begin{equation*}
\left\{\begin{aligned}
&
\frac{E_{3,j-1}^{(m,n)}-2E_{3,j}^{(m,n)}+E_{3,j+1}^{(m,n)}}{h^2}+\Big{(}\kappa_0^2-\big{(}\frac{m\pi}{a}\big{)}^2-\big{(}\frac{n\pi}{b}\big{)}^2\Big{)}E_{3,j}^{(m,n)}=0,~ j=1,2,\cdots,J,\\
&E_{3,1}^{(m,n)}=E_{3,0}^{(m,n)}.
\end{aligned}\right.
\end{equation*}
The above discrete systems can be written in the matrix form
\begin{equation}\label{dfE1M}
\big{(}\boldsymbol{A}_1+\boldsymbol{D}^{(m,n)}\big{)}\boldsymbol{E}_l^{(m,n)}+\boldsymbol{a}_{J}E_{l,J+1}^{(m,n)}=0,~(m,n)\in\mathbb{N}^2_l, l = 1,2,
\end{equation}
and
\begin{equation}\label{dfE3M}
\big{(}\boldsymbol{A}_2+\boldsymbol{D}^{(m,n)}\big{)}\boldsymbol{E}_3^{(m,n)}+\boldsymbol{a}_{J}E_{3,J+1}^{(m,n)}=0,~(m,n)\in\mathbb{N}^2_3,
\end{equation}
where the vectors of unknowns $
\boldsymbol{E}_l^{(m,n)}=\Big{(}E_{l,1}^{(m,n)},E_{l,2}^{(m,n)},\cdots,E_{l,J}^{
(m,n)}\Big{)}^\top,\quad l=1,2,3$,
\begin{equation*}
  \boldsymbol{A}_1=\begin{pmatrix}
  -2&1&&\\
  1&-2&1&\\
  &\ddots&\ddots&\ddots\\
  &&1&-2
  \end{pmatrix},\quad
  \boldsymbol{A}_2=\begin{pmatrix}
  -1&1&&\\
  1&-2&1&\\
  &\ddots&\ddots&\ddots\\
  &&1&-2
  \end{pmatrix},\quad
  \boldsymbol{a}_{J}=\begin{pmatrix}
  0\\
  \vdots\\
  0\\
  1
  \end{pmatrix},
\end{equation*}
and
\begin{equation*}
\boldsymbol{D}^{(m,n)}=h^2\Big{(}\kappa_0^2-\big{(}\frac{m\pi}{a}\big{)}^2-\big{(}\frac{n\pi}{b}\big{)}^2\Big{)}\boldsymbol{I}_J,
\end{equation*}
Here $\boldsymbol{I}_J$  is the $J\times J$ identity matrix.\par

Next, we discuss the discretization of the transparent boundary condition
\eqref{ENU1}. A simple calculation from the first component of \eqref{ENU1}
yields
\begin{align}
\nonumber   \frac{\partial E_3}{\partial x_1}-\frac{\partial E_1}{\partial
x_3}&=2(i\alpha_1p_3+i\beta p_1)e^{i(\alpha_1x_1+\alpha_2x_2)}
+2\kappa_0^2\int_{\Gamma}E_1(\boldsymbol{y})g(\boldsymbol{x},\boldsymbol{y}){
\rm d}s_{\boldsymbol{y}} \\
\label{ENUC1}
&\quad +2\int_{\Gamma}\big{(}-\partial_{y_1}E_2(\boldsymbol{y})+\partial_{y_2}
E_1(\boldsymbol{y})\big{)}\partial_{x_2}g(\boldsymbol{x},\boldsymbol{y}){\rm
d}s_{\boldsymbol{y}}.
\end{align}
Substituting \eqref{E123F} into \eqref{ENUC1}, we have
\begin{align}
\nonumber   &\sum\limits_{k\in\mathbb{N}_3^2}E_3^{(k)}(0)\frac{k_1\pi}{a}
\cos\Big{(}\frac{k_1\pi x_1}{a}\Big{)}\sin\Big{(}\frac{k_2\pi
x_2}{b}\Big{)}-\sum\limits_{k\in\mathbb{N}_1^2}\frac{{\rm d}
E_1^{(k)}(0)}{{\rm d} x_3}\cos\Big{(}\frac{k_1\pi x_1}{a}\Big{)}
\sin\Big{(}\frac{k_2\pi x_2}{b}\Big{)} \\
 \nonumber  &=2(i\alpha_1p_3+i\beta p_1)e^{i(\alpha_1x_1+\alpha_2x_2)} \\
\nonumber
&\quad +2\kappa_0^2\sum\limits_{k\in\mathbb{N}_1^2}E_1^{(k)}(0)\int_{\Gamma}
\cos\Big{(}\frac{k_1\pi y_1}{a}\Big{)}\sin\Big{(}\frac{k_2\pi
y_2}{b}\Big{)}g(\boldsymbol{x},\boldsymbol{y}){\rm d}s_{\boldsymbol{y}}\\
\nonumber&\quad -2\sum\limits_{k\in\mathbb{N}_2^2}E_2^{(k)}(0)\frac{k_1\pi}{a}
\int_{\Gamma}\cos\Big{(}\frac{k_1\pi y_1}{a}\Big{)}\cos\Big{(}\frac{k_2\pi
y_2}{b}\Big{)}\partial_{x_2}g(\boldsymbol{x},\boldsymbol{y}){\rm
d}s_{\boldsymbol{y}}\\\label{ENUC11}
&\quad
+2\sum\limits_{k\in\mathbb{N}_1^2}E_1^{(k)}(0)\frac{k_2\pi}{b}\int_{\Gamma }
\cos\Big{(}\frac{k_1\pi y_1}{a}\Big{)}\cos\Big{(}\frac{k_2\pi
y_2}{b}\Big{)}\partial_{x_2}g(\boldsymbol{x},\boldsymbol{y}){\rm
d}s_{\boldsymbol{y}}.
\end{align}
Multiplying both sides of \eqref{ENUC11} by ${\rm cos}\big{(}\frac{m\pi x_1}{a}\big{)}{\rm sin}\big{(}\frac{n\pi x_2}{b}\big{)}, (m,n)\in\mathbb{N}^2_1$ and integrating over $\Gamma$, we obtain
\begin{align*}
E_3^{(m,n)}(0)p^{(m,n)}-\frac{{\rm d}
E_1^{(m,n)}(0)}{{\rm d} x_3}q^{(m,n)}
 =2(i\alpha_1p_3+i\beta p_1)\tilde{g}_1^{(m,n)}
+2\kappa_0^2\sum\limits_{k\in\mathbb{N}_1^2}E_1^{(k)}(0)\tilde{F}_{1,(k)}^{(m,n)
}\\
-2\sum\limits_{k\in\mathbb{N}_2^2}E_2^{(k)}(0)\frac{k_1\pi}{a}\tilde{G}_{1,
(k)}^{(m,n)}
+2\sum\limits_{k\in\mathbb{N}_1^2}E_1^{(k)}(0)\frac{k_2\pi}{b}\tilde{H}_{1,(k)}^
{(m,n)},
\end{align*}
where
\begin{align*}
p^{(m,n)} = \begin{cases}
0, &\mbox{ if } m=0,\\
\frac{bm\pi}{4}, &\mbox{ others, }
\end{cases}
\quad
q^{(m,n)} = \begin{cases}
\frac{ab}{2}, &\mbox{ if } m=0,\\
\frac{ab}{4}, &\mbox{ others, }
\end{cases}
\end{align*}
and
\begin{align}
&\tilde{g}_1^{(m,n)}=\int_{\Gamma}\cos\Big{(}\frac{m\pi x_1}{a}\Big{)}
\sin\Big{(}\frac{n\pi x_2}{b}\Big{)}e^{i(\alpha_1x_1+\alpha_2x_2)}{\rm
d}s_{\boldsymbol{x}},\label{nosig}\\
& \tilde{F}_{1,(k)}^{(m,n)}=\int_{\Gamma}\cos\Big{(}\frac{m\pi
x_1}{a}\Big{)}\sin\Big{(}\frac{n\pi x_2}{b}\Big{)}\Big{(}\int_{\Gamma}
\cos\Big{(}\frac{k_1\pi y_1}{a}\Big{)}\sin\Big{(}\frac{k_2\pi
y_2}{b}\Big{)}g(\boldsymbol{x},\boldsymbol{y}){\rm
d}s_{\boldsymbol{y}}\Big{)}{\rm d}s_{\boldsymbol{x}},\label{sig1}\\
&\tilde{G}_{1,(k)}^{(m,n)}=\int_{\Gamma}\cos\Big{(}\frac{m\pi
x_1}{a}\Big{)}\sin\Big{(}\frac{n\pi x_2}{b}\Big{)}\Big{(}\int_{\Gamma}
\cos\Big{(}\frac{k_1\pi y_1}{a}\Big{)}\cos\Big{(}\frac{k_2\pi
y_2}{b}\Big{)}\partial_{x_2}g(\boldsymbol{x},\boldsymbol{y}){\rm
d}s_{\boldsymbol{y}}\Big{)}{\rm d}s_{\boldsymbol{x}},\label{sig2}\\
&\tilde{H}_{1,(k)}^{(m,n)}=\int_{\Gamma}\cos\Big{(}\frac{m\pi
x_1}{a}\Big{)}\sin\Big{(}\frac{n\pi x_2}{b}\Big{)}\Big{(}\int_{\Gamma}
\cos\Big{(}\frac{k_1\pi y_1}{a}\Big{)}\cos\Big{(}\frac{k_2\pi
y_2}{b}\Big{)}\partial_{x_2}g(\boldsymbol{x},\boldsymbol{y}){\rm
d}s_{\boldsymbol{y}}\Big{)}{\rm d}s_{\boldsymbol{x}}. \label{sig3}
\end{align}

By using a backward finite difference scheme for the normal derivative and the fact that $E_{l,J+1}^{(k)}=E_l^{(k)}(0)$, $l=1,2,3$, we get
\begin{align} \label{c1mn}
E_{3,J+1}^{(m,n)}p^{(m,n)} -\frac{E_{1,J+1}^{(m,n)}-E_{1,J}^{
(m,n)}}{h}q^{(m,n)} =2(i\alpha_1p_3+i\beta p_1)\tilde{g}_1^{(m,n)}
+2\kappa_0^2\sum\limits_{k\in\mathbb{N}^2_1}E_{1,J+1}^{(k)}\tilde{F}_{1,(k)}^{(m
,n)} \notag\\
-2\sum\limits_{k\in\mathbb{N}^2_2}E_{2,J+1}^{(k)}\Big(\frac{k_1\pi}{a}
\Big)\tilde{G}_{1,(k)}^{(m,n)}
+2\sum\limits_{k\in\mathbb{N}^2_1}E_{1,J+1}^{(k)}\Big(\frac{k_2\pi}{b}
\Big)\tilde{H}_{1, (k)}^{(m,n)}.
\end{align}
For $(m,n)\in\mathbb{N}^2_1$, we define the following notations:
\begin{align*}
&g_1^{(m,n)}:=\frac{h}{q^{(m,n)}}2(i\alpha_1p_3+i\beta p_1)\tilde{g}_1^{(m,n)},\\
& F_{1,(k)}^{(m,n)}:=\frac{h}{q^{(m,n)}}2\kappa_0^2\tilde{F}_{1,(k)}^{(m,n)},\\
&G_{1,(k)}^{(m,n)}:=\frac{h}{q^{(m,n)}}\frac{-2k_1\pi}{a}\tilde{G}_{1,(k)}^{(m,n)},\\
&H_{1,(k)}^{(m,n)}:=\frac{h}{q^{(m,n)}}\frac{2k_2\pi}{b}\tilde{H}_{1,(k)}^{(m,n)}.
\end{align*}
Thus, we obtain from \eqref{c1mn} that
\begin{align*}
\nonumber  &
E_{3,J+1}^{(m,n)}\frac{p^{(m,n)}h}{q^{(m,n)}}-E_{1,J+1}^{(m,n)}+E_{1,J}^{(m,n)}
-\sum\limits_{k\in\mathbb{N}^2_1}E_{1,J+1
}^{(k)}F_{1,(k)}^{(m,n)}  \\
&\quad
-\sum\limits_{k\in\mathbb{N}^2_2}E_{2,J+1}^{(k)}G_{1,(k)}^{(m,n)}
-\sum\limits_ {
k\in\mathbb{N}^2_1}E_{1,J+1}^{(k)}H_{1,(k)}^{(m,n)}=g_1^{(m,n)},\quad
(m,n)\in\mathbb{N}^2_1,
\end{align*}
which can be written in a matrix form
\begin{equation}\label{ENU14}
\hat{\boldsymbol{I}}_1\boldsymbol{E}_{1,J}+(-\hat{\boldsymbol{I}}_1-\boldsymbol{
F}_1-\boldsymbol{H}_1)\boldsymbol{E}_{1,J+1}-\boldsymbol{G}_1\boldsymbol{E}_{2,
J+1}+\boldsymbol{I}_1\boldsymbol{E}_{3,J+1}=\boldsymbol{g}_{1}.
\end{equation}
Here
$\hat{\boldsymbol{I}}_1=\boldsymbol{I}_{((M+1)N)}$, $\boldsymbol{I}_1=\tilde{
\boldsymbol{I}}_1\otimes\boldsymbol{I}_N$, $\otimes$ denotes the
Kronecker product, and
\begin{equation*}
  \tilde{\boldsymbol{I}}_1=\begin{pmatrix}
  0&\cdots&0\\
  \frac{\pi h}{a}&&\\
  &\ddots&\\
  &&\frac{M\pi h}{a}
  \end{pmatrix}.
\end{equation*}
For clarity, we refer to Appendix A for the entries of
$\boldsymbol{F}_1,~\boldsymbol{H}_1,\boldsymbol{G}_1$,
$\boldsymbol{g}_{1}$, and $\boldsymbol{E}_{l,j}$ for $l=1,2,3, 0\le j\le
J+1$.

Similarly, the second component of TBC \eqref{ENU1} can be discretized as
\begin{equation}\label{ENU24}
\hat{\boldsymbol{I}}_2\boldsymbol{E}_{2,J}-\boldsymbol{H}_2\boldsymbol{E}_{1,J+1}+(-\hat{\boldsymbol{I}}_2-\boldsymbol{F}_2-\boldsymbol{G}_2)\boldsymbol{E}_{2,J+1}+\boldsymbol{I}_{2}\boldsymbol{E}_{3,J+1}=\boldsymbol{g}_{2},
\end{equation}
where $\hat{\boldsymbol{I}}_2=\boldsymbol{I}_{(M(N+1))}$,
$\boldsymbol{I}_2=\boldsymbol{I}_M\otimes\tilde{\boldsymbol{I}}_2$, and
\begin{equation*}
  \tilde{\boldsymbol{I}}_2=\begin{pmatrix}
  0&\cdots&0\\
  \frac{\pi h}{b}&&\\
  &\ddots&\\
  &&\frac{N\pi h}{b}
  \end{pmatrix}_{(N+1)\times N}.
\end{equation*}
Again, the entries of the vectors
$\boldsymbol{F}_2,~\boldsymbol{H}_2,\boldsymbol{G}_2$ and
$\boldsymbol{g}_{2}$ can be found in Appendix A.

Recall the divergence free condition on the surface $\Gamma$,
\begin{equation}\label{ENUC3}
 \partial_{x_1}E_1+\partial_{x_2}E_2+\partial_{x_3}E_3=0.
\end{equation}
Substituting \eqref{E123F} into \eqref{ENUC3}, we have
\begin{align}
 \nonumber
\sum\limits_{k\in\mathbb{N}_1^2}E_1^{(k)}(0)\Big(\frac{-k_1\pi}{a}\Big)
\sin\Big{(}\frac{k_1\pi x_1}{a}\Big{)}\sin\Big{(}\frac{k_2\pi
x_2}{b}\Big{)}
+\sum\limits_{k\in\mathbb{N}_2^2}E_2^{(k)}(0)\Big(\frac{-k_2\pi}{b }\Big)
\sin\Big{(}\frac{k_1\pi x_1}{a}\Big{)}\sin\Big{(}\frac{k_2\pi
x_2}{b}\Big{)}\\\label{ENUC31}
+\sum\limits_{k\in\mathbb{N}_3^2}\frac{\partial E_3^{(k)}(0)}{\partial x_3}
\sin\Big{(}\frac{k_1\pi x_1}{a}\Big{)}\sin\Big{(}\frac{k_2\pi
x_2}{b}\Big{)}=0.
\end{align}
Multiplying both side of \eqref{ENUC31} by ${\rm sin}\big{(}\frac{m\pi
x_1}{a}\big{)}{\rm sin}\big{(}\frac{n\pi x_2}{b}\big{)}, (m,n)\in
\mathbb{N}^2_3$, integrating over $\Gamma$, and using the orthogonality of the
trigonometric functions, we obtain
\begin{equation}\label{ENU32}
\left(\frac{-m\pi}{a}\frac{ab}{4}\right)
E_{1}^{(m,n)}(0)+\left(\frac{-n\pi}{b}\frac{ab}{4}\right)E_{2}^{(m,n)}
(0)+\left(\frac{ab}{4}\right)\frac{\partial E_3^{(m,n)}(0)}{\partial x_3}=0.
\end{equation}
By using a backward finite difference scheme, we get
\begin{equation}\label{ENU33}
E_{3,J}^{(m,n)}+\left(\frac{m\pi}{a}\right)E_{1,J+1}^{(m,n)}+\left(\frac{n\pi}{b
}\right)E_{2,J+1}^{(m,n)} - E_{3,J+1}^{(m,n)}=0.
\end{equation}
Let
\begin{equation*}
  \tilde{\boldsymbol{I}}_3=\begin{pmatrix}
  0&\frac{\pi h}{a}&&\\
  \vdots&&\ddots&\\
  0&&&\frac{M\pi h}{a}
  \end{pmatrix},\quad
  \tilde{\boldsymbol{I}}_4=\begin{pmatrix}
  0&\frac{\pi h}{b}&&\\
  \vdots&&\ddots&\\
  0&&&\frac{N\pi h}{b}
  \end{pmatrix},
\end{equation*}
$\boldsymbol{F}_3=\tilde{\boldsymbol{I}}_3\otimes\boldsymbol{I}_N$ and
$\boldsymbol{G}_3=\boldsymbol{I}_M\otimes\tilde{\boldsymbol{I}}_4$.
The discrete system \eqref{ENU33} can be rewritten as
\begin{equation}\label{ENU34}
\hat{\boldsymbol{I}}_3\boldsymbol{E}_{3,J}+\boldsymbol{F}_3\boldsymbol{E}_{1,J+1}+\boldsymbol{G}_3\boldsymbol{E}_{2,J+1}-\hat{\boldsymbol{I}}_3\boldsymbol{E}_{3,J+1}=0,
\end{equation}
where $\hat{\boldsymbol{I}}_3=\boldsymbol{I}_{(MN)}$.
It follows from \eqref{ENU14}--\eqref{ENU24} and \eqref{ENU34} that
\begin{equation}\label{ENUGa}
  \begin{pmatrix}
  \hat{\boldsymbol{I}}_1&&\\
  &\hat{\boldsymbol{I}}_2&\\
  &&\hat{\boldsymbol{I}}_3
  \end{pmatrix}
    \begin{pmatrix}
  \boldsymbol{E}_{1,J}\\
  \boldsymbol{E}_{2,J}\\
  \boldsymbol{E}_{3,J}
  \end{pmatrix}
  +
  \begin{pmatrix}
  -\hat{\boldsymbol{I}}_1-\boldsymbol{F}_1-\boldsymbol{H}_1&-\boldsymbol{G}_1&\boldsymbol{I}_1\\
  -\boldsymbol{H}_2&-\hat{\boldsymbol{I}}_2-\boldsymbol{F}_2-\boldsymbol{G}_2&\boldsymbol{I}_2\\
  \boldsymbol{F}_3&\boldsymbol{G}_3&-\hat{\boldsymbol{I}}_3
  \end{pmatrix}
    \begin{pmatrix}
  \boldsymbol{E}_{1,J+1}\\
  \boldsymbol{E}_{2,J+1}\\
  \boldsymbol{E}_{3,J+1}
  \end{pmatrix}
  =
  \begin{pmatrix}
  \boldsymbol{g}_1\\
  \boldsymbol{g}_2\\
  0
  \end{pmatrix}
\end{equation}

Clearly, the linear systems \eqref{dfE1M}--\eqref{dfE3M} and \eqref{ENUGa} are
coupled and give the global system. Next, we use Gaussian elimination
method to decouple the global system into a linear system with the unknowns only
on the aperture, which may reduce the computational complexity greatly and lead
to a fast algorithm.

Let
\begin{equation}\label{LU1}
  \boldsymbol{L}_1^{(m,n)}\boldsymbol{U}_1^{(m,n)}=\boldsymbol{A}_1+\boldsymbol{D}^{(m,n)},\quad (m,n)\in\mathbb{N}^2_l, ~l = 1,2,
\end{equation}
and
\begin{equation}\label{LU2}
  \boldsymbol{L}_2^{(m,n)}\boldsymbol{U}_2^{(m,n)}=\boldsymbol{A}_2+\boldsymbol{D}^{(m,n)},\quad (m,n)\in\mathbb{N}^2_3,
\end{equation}
be the LU-decomposition, where $\boldsymbol{A}_1+\boldsymbol{D}^{(m,n)}$ and
$\boldsymbol{A}_2+\boldsymbol{D}^{(m,n)}$ are the symmetric tridiagonal matrices
in \eqref{dfE1M} and \eqref{dfE3M}, respectively. Since
$\boldsymbol{L}_1^{(m,n)}$ and $\boldsymbol{L}_2^{(m,n)}$ are nonsingular, we
obtain
\begin{equation}\label{LUE1}
\boldsymbol{U}_1^{(m,n)}\boldsymbol{E}_l^{(m,n)}+\big{(}\boldsymbol{L}_1^{(m,n)}\big{)}^{-1}\boldsymbol{a}_{J+1}E_{l,J+1}^{(m,n)}=0,\quad (m,n)\in\mathbb{N}^2_l, ~l = 1,2,
\end{equation}
\begin{equation}\label{LUE3}
\boldsymbol{U}_2^{(m,n)}\boldsymbol{E}_3^{(m,n)}+\big{(}\boldsymbol{L}_2^{(m,n)}\big{)}^{-1}\boldsymbol{a}_{J+1}E_{3,J+1}^{(m,n)}=0,\quad (m,n)\in\mathbb{N}^2_3,
\end{equation}
where $\boldsymbol{U}_1^{(m,n)}=\big{(}r_{1,(pq)}^{m,n}\big{)}$ and $\boldsymbol{U}_2^{(m,n)}=\big{(}r_{2,(pq)}^{m,n}\big{)}$.

Combining the last equations of the systems \eqref{LUE1} and \eqref{LUE3} gives
\begin{equation}\label{EJJ1}
  \begin{pmatrix}
  \boldsymbol{R}_1&&\\
  &\boldsymbol{R}_2&\\
  &&\boldsymbol{R}_3
  \end{pmatrix}
    \begin{pmatrix}
  \boldsymbol{E}_{1,J}\\
  \boldsymbol{E}_{2,J}\\
  \boldsymbol{E}_{3,J}
  \end{pmatrix}
  +
  \begin{pmatrix}
  \hat{\boldsymbol{I}}_1&&\\
  &\hat{\boldsymbol{I}}_2&\\
  &&\hat{\boldsymbol{I}}_3
  \end{pmatrix}
    \begin{pmatrix}
  \boldsymbol{E}_{1,J+1}\\
  \boldsymbol{E}_{2,J+1}\\
  \boldsymbol{E}_{3,J+1}
  \end{pmatrix}
  =
 0,
\end{equation}
where
\begin{equation*}
  \boldsymbol{R}_l={\rm diag}\big{(}r_{1,(JJ)}^{(m,n)}\big{)}, \quad(m,n)\in\mathbb{N}^2_l,~ l = 1,2,
\end{equation*}
\begin{equation*}
  \boldsymbol{R}_3={\rm diag}\big{(}r_{2,(JJ)}^{(m,n)}\big{)}, \quad(m,n)\in\mathbb{N}^2_3.
\end{equation*}
If $\kappa_0^2$ is not an eigenvalue of the Helmholtz operator with Dirichlet
boundary condition, the continuous Helmholtz problem admits a unique
solution; for $h$ small enough, as an approximate problem, the discrete
Helmholtz problem can also be shown to have a unique solution \cite{YGL20},
which implies that
\begin{equation}\label{r1}
r_{1,(JJ)}^{(m,n)}\neq 0, \quad (m,n)\in\mathbb{N}^2_l, ~l = 1,2,
\end{equation}
and
\begin{equation}\label{r2}
r_{2,(JJ)}^{(m,n)}\neq 0, \quad (m,n)\in\mathbb{N}^2_3.
\end{equation}
Consequently, combining \eqref{EJJ1} and \eqref{ENUGa} yields
\begin{equation}\label{EJ1}
  \begin{pmatrix}
  -\hat{\boldsymbol{I}}_1-\boldsymbol{F}_1-\boldsymbol{H}_1-\boldsymbol{R}_1^{-1}&-\boldsymbol{G}_1&\boldsymbol{I}_1\\
  -\boldsymbol{H}_2&-\hat{\boldsymbol{I}}_2-\boldsymbol{F}_2-\boldsymbol{G}_2-\boldsymbol{R}_2^{-1}&\boldsymbol{I}_2\\
  \boldsymbol{F}_3&\boldsymbol{G}_3&-\hat{\boldsymbol{I}}_3-\boldsymbol{R}_3^{-1}
  \end{pmatrix}
    \begin{pmatrix}
  \boldsymbol{E}_{1,J+1}\\
  \boldsymbol{E}_{2,J+1}\\
  \boldsymbol{E}_{3,J+1}
  \end{pmatrix}
  =
  \begin{pmatrix}
  \boldsymbol{g}_1\\
  \boldsymbol{g}_2\\
  0
  \end{pmatrix}.
\end{equation}
Solving the linear system \eqref{EJ1} gives the solution $E_{l,J+1},
l=1, 2, 3$ on the interface $\Gamma$. The rest of the
unknowns can be simply obtained by solving the following systems:
\begin{equation}\label{EInter}
  \begin{split}
&\big{(}\boldsymbol{A}_1+\boldsymbol{D}^{(m,n)}\big{)}\boldsymbol{E}_l^{(m,n)}
=-\boldsymbol{a}_{J+1}E_{l,J+1}^{(m,n)}, \quad l = 1,2,\\
  &\big{(}\boldsymbol{A}_2+\boldsymbol{D}^{(m,n)}\big{)}\boldsymbol{E}_3^{(m,n)}=-\boldsymbol{a}_{J+1}E_{3,J+1}^{(m,n)}.
  \end{split}
\end{equation}

\begin{remark}
Since the medium is assumed to be homogeneous in the cavity, it follows from the
Maxwell equation \eqref{beq1} that the electrical field $\boldsymbol{E}$ is
divergence free in $D$. Although the solutions are solved separately in $D$,
they admit the series expansions \eqref{E123F} and satisfy the divergence free
condition due to \eqref{beq1}.
\end{remark}

\section{Layered media}

\begin{figure}
\includegraphics[width=0.5\textwidth]{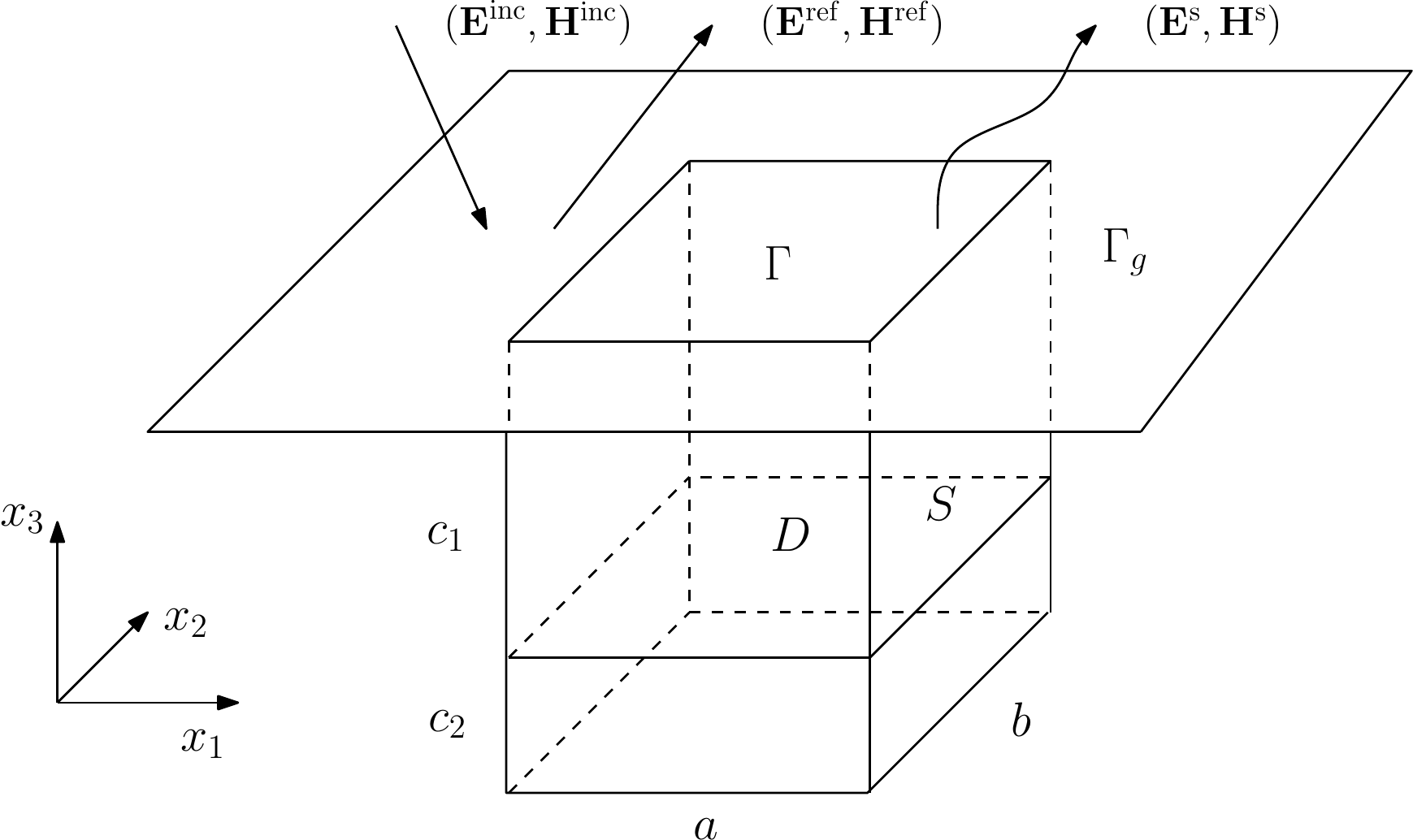}
\caption{The problem geometry of the electromagnetic scattering by a
rectangular cavity filled with a layered medium.}
\label{fig:2}
\end{figure}

This section is devoted to the numerical solution of the electromagnetic
scattering by an open cavity with a layered medium. Specifically, we assume that
the cavity is filled with a multi-layered medium, which is characterized by the
piecewise constant dielectric permittivity $\varepsilon_l, l=1,2,\cdots, L$. The
medium is still assumed to be nonmagnetic with a constant magnetic permeability
$\mu=\mu_0$ everywhere and has a constant dielectric permittivity
$\varepsilon=\varepsilon_0$ in the upper half space. Without loss of generality,
we discuss a two-layered medium in $D$. Denote by $c_1$ and $c_2$ the depth of
the two layer domain $D_1$ and $D_2$, respectively. The problem geometry is
depicted in Figure \ref{fig:2}. The open aperture of the cavity
$\Gamma=[0,a]\times[0,b]$ and the total depth of the cavity is $c$, i.e.,
$c=c_1+c_2$.

Let $\boldsymbol{E}_1=(u_1,u_2,u_3)$ and $\boldsymbol{E}_2=(v_1,v_2,v_3)$ be
the total electric field in domain $D_1$ and $D_2$, respectively.
Similar to the homogeneous case, it can be shown from the boundary condition and
divergence free condition that $u_j$ and $v_j, j=1,2,3$ admit the following
Fourier series expansions:
\begin{equation}\label{E1F}
  \left\{\begin{aligned}
  u_1(x_1,x_2,x_3)&=\sum\limits_{k\in
\mathbb{N}^2}u_1^{(k)}(x_3)\cos\Big{(}\frac{k_1\pi
x_1}{a}\Big{)}\sin\Big{(}\frac{k_2\pi x_2}{b}\Big{)},\\
  u_2(x_1,x_2,x_3)&=\sum\limits_{k\in
\mathbb{N}^2}u_2^{(k)}(x_3)\sin\Big{(}\frac{k_1\pi x_1
}{a}\Big{)}\cos\Big{(}\frac{k_2\pi x_2}{b}\Big{)},\\
  u_3(x_1,x_2,x_3)&=\sum\limits_{k\in
\mathbb{N}^2}u_3^{(k)}(x_3)\sin\Big{(}\frac{k_1\pi x_1
}{a}\Big{)}\sin\Big{(}\frac{k_2\pi x_2}{b}\Big{)},
  \end{aligned}\right.
  \end{equation}
and
\begin{equation}\label{E2F}
 \left\{ \begin{aligned}
  v_1(x_1,x_2,x_3)&=\sum\limits_{k\in
\mathbb{N}^2}v_1^{(k)}(x_3)\cos\Big{(}\frac{k_1\pi
x_1}{a}\Big{)}\sin\Big{(}\frac{k_2\pi x_2}{b}\Big{)},\\
  v_2(x_1,x_2,x_3)&=\sum\limits_{k\in
\mathbb{N}^2}v_2^{(k)}(x_3)\sin\Big{(}\frac{k_1\pi x_1
}{a}\Big{)}\cos\Big{(}\frac{k_2\pi x_2}{b}\Big{)},\\
  v_3(x_1,x_2,x_3)&=\sum\limits_{k\in
\mathbb{N}^2}v_3^{(k)}(x_3)\sin\Big{(}\frac{k_1\pi x_1
}{a}\Big{)}\sin\Big{(}\frac{k_2\pi x_2}{b}\Big{)},
  \end{aligned}\right.
  \end{equation}
where $k=(k_1,k_2)\in \mathbb{N}^2$.

In the lower part of the layered medium $D_2$, the electric field
$\boldsymbol{E}_2=(v_1,v_2,v_3)$ satisfies the Helmholtz equation
\begin{equation}\label{E2Heq}
  \Delta \boldsymbol{E}_2+\kappa_2^2\boldsymbol{E}_2=0 \quad \mathrm{in}~D_2,
\end{equation}
the homogeneous Dirichlet boundary condition
\begin{equation}\label{v12D}
  v_1(x_1,x_2,-c)=v_2(x_1,x_2,-c)=0,
\end{equation}
and the homogeneous Neumann boundary condition
\begin{equation}\label{v3N}
\partial_{x_3}v_3(x_1,x_2,-c)=0,
\end{equation}
where $\kappa_2=\omega\sqrt{\varepsilon_2\mu}$ is the wavenumber in $D_2$.

Substituting \eqref{E2F} into \eqref{E2Heq}--\eqref{v3N}, we get the
second order ordinary differential equations with the homogeneous Dirichlet
boundary condition at $x_3=-c$ for the Fourier coefficients $
v_l^{(m,n)}, l=1, 2$:
\begin{equation}\label{odev1}
\left\{\begin{aligned}
 & \frac{{\rm d}^2} {{\rm d}x_3^2} v_l^{(m,n)}(x_3)+\Big{(}\kappa_2^2-\big{(}\frac{m\pi}{a}\big{)}^2-\big{(}\frac{n\pi}{b}\big{)}^2\Big{)}v_l^{(m,n)}(x_3)=0,\quad x_3\in(-c,-c_1),\\
 &v_l^{(m,n)}(-c)=0,
 \end{aligned}\right.
\end{equation}
where $(m,n)\in\mathbb{N}^2_l$, and the second order ordinary differential
equations with the homogeneous Neumann boundary condition at $x_3=-c$ for the
Fourier coefficients $v_3^{(m,n)}$:
\begin{equation}\label{odev3}
\left\{\begin{aligned}
&\frac{{\rm d}^2}{{\rm d}x_3^2}v_3^{(m,n)}(x_3)+\Big{(}\kappa_2^2-\big{(}\frac{m\pi}{a}\big{)}^2-\big{(}\frac{n\pi}{b}\big{)}^2\Big{)}v_3^{(m,n)}(x_3)=0,\quad x_3\in(-c,-c_1),\\
&\frac{\rm d}{{\rm d}x_3}v_3^{(m,n)}(-c)=0.
   \end{aligned}\right.
\end{equation}
where $(m,n)\in\mathbb{N}^2_3$.

Define by $\{x_3^j\}_{j=0}^{j=J+1}$ a set of uniformly distributed grid points
in $[-c,-c_1]$, where $h=x_3^{j+1}-x_3^j$. Let $v_{l,j}^{(m,n)}, l=1,2,3 $ be
the finite difference solution of $v_l^{(m,n)}(x_3)$ at the point $x_3=x_3^j$.
Similar to the discretization of \eqref{ode1x3}--\eqref{ode3x3},
the discrete system of \eqref{odev1}--\eqref{odev3} can be written in the matrix
form
\begin{equation}\label{dfv1M}
\big{(}\boldsymbol{A}_1+\boldsymbol{D}_2^{(m,n)}\big{)}\boldsymbol{v}_l^{(m,n)}+\boldsymbol{a}_{J}v_{l,J+1}^{(m,n)}=0,~(m,n)\in\mathbb{N}^2_l, l = 1,2,
\end{equation}
and
\begin{equation}\label{dfv3M}
\big{(}\boldsymbol{A}_2+\boldsymbol{D}_2^{(m,n)}\big{)}\boldsymbol{v}_3^{(m,n)}+\boldsymbol{a}_{J}v_{3,J+1}^{(m,n)}=0,~(m,n)\in\mathbb{N}^2_3,
\end{equation}
where the vector of unknowns $
\boldsymbol{v}_l^{(m,n)}=\Big{(}v_{l,1}^{(m,n)},v_{l,2}^{(m,n)},\cdots,v_{l,J}^{
(m,n)}\Big{)}^\top,\quad l=1,2,3$,
\begin{equation*}
  \boldsymbol{A}_1=\begin{pmatrix}
  -2&1&&\\
  1&-2&1&\\
  &\ddots&\ddots&\ddots\\
  &&1&-2
  \end{pmatrix},\quad
  \boldsymbol{A}_2=\begin{pmatrix}
  -1&1&&\\
  1&-2&1&\\
  &\ddots&\ddots&\ddots\\
  &&1&-2
  \end{pmatrix},\quad
  \boldsymbol{a}_{J}=\begin{pmatrix}
  0\\
  \vdots\\
    0\\
  1
  \end{pmatrix},
\end{equation*}
and
\begin{equation*}
\boldsymbol{D}_2^{(m,n)}=h^2\Big{(}\kappa_2^2-\big{(}\frac{m\pi}{a}\big{)}^2-\big{(}\frac{n\pi}{b}\big{)}^2\Big{)}\boldsymbol{I}_J,
\end{equation*}

Again, we apply the Gaussian elimination method to solve the linear system
\eqref{dfv1M}--\eqref{dfv3M}. Let
\begin{equation}\label{LUd21}
  \boldsymbol{L}_1^{(m,n)}\boldsymbol{U}_1^{(m,n)}=\boldsymbol{A}_1+\boldsymbol{D}_2^{(m,n)},\quad (m,n)\in\mathbb{N}^2_l, ~ l =1,2,
\end{equation}
and
\begin{equation}\label{LUd22}
  \boldsymbol{L}_2^{(m,n)}\boldsymbol{U}_2^{(m,n)}=\boldsymbol{A}_2+\boldsymbol{D}_2^{(m,n)},\quad (m,n)\in\mathbb{N}^2_3,
\end{equation}
be the LU-decomposition, where
$\boldsymbol{U}_1^{(m,n)}=\big{(}r_{1,(pq)}^{(m,n)}\big{)}$ and
$\boldsymbol{U}_2^{(m,n)}=\big{(}r_{2,(pq)}^{(m,n)}\big{)}$. Since
$\boldsymbol{L}_1^{(m,n)}$ and $\boldsymbol{L}_2^{(m,n)}$ are nonsingular, we
obtain
\begin{equation}\label{LUv1}
\boldsymbol{U}_1^{(m,n)}\boldsymbol{v}_l^{(m,n)}+\big{(}\boldsymbol{L}_1^{(m,n)}\big{)}^{-1}\boldsymbol{a}_{J}v_{l,J+1}^{(m,n)}=0,\quad (m,n)\in\mathbb{N}^2_l,~ l = 1,2,
\end{equation}
and
\begin{equation}\label{LUv3}
\boldsymbol{U}_2^{(m,n)}\boldsymbol{v}_3^{(m,n)}+\big{(}\boldsymbol{L}_2^{(m,n)}\big{)}^{-1}\boldsymbol{a}_{J}v_{3,J+1}^{(m,n)}=0,\quad (m,n)\in\mathbb{N}^2_3.
\end{equation}
Combining the last equations of the systems \eqref{LUv1} and \eqref{LUv3} gives
\begin{align}
\label{v1JJ1}&r_{1,(JJ)}^{(m,n)}v_{l,J}^{(m,n)}+v_{l,J+1}^{(m,n)}=0,\quad
(m,n)\in\mathbb{N}^2_l, ~l =1,2,\\
\label{v3JJ1}&r_{2,(JJ)}^{(m,n)}v_{3,J}^{(m,n)}+v_{3,J+1}^{(m,n)}=0,\quad (m,n)\in\mathbb{N}^2_3.
\end{align}

In the upper part of the layered medium $D_1$, the electric field
$\boldsymbol{E}_1=(u_1,u_2,u_3)$
satisfies the Helmholtz equation
\begin{equation}\label{E1Heq}
  \Delta \boldsymbol{E}_1+\kappa_1^2\boldsymbol{E}_1=0 \quad \mathrm{in}~D_1,
\end{equation}
where $\kappa_1=\omega\sqrt{\varepsilon_1\mu}$ is the wavenumber in $D_1$.
Substituting \eqref{E1F} into \eqref{E1Heq} yields
\begin{equation}\label{odeu1}
\frac{{\rm d}^2} {{\rm d}x_3^2} u_l^{(m,n)}(x_3)+\Big{(}\kappa_1^2-\big{(}\frac{m\pi}{a}\big{)}^2-\big{(}\frac{n\pi}{b}\big{)}^2\Big{)}u_l^{(m,n)}(x_3)=0,\quad x_3\in(-c_1,0),\quad (m,n)\in\mathbb{N}_l^2,
\end{equation}
for  $l =1,2,3.$

Let $\{x_3^i\}_{i=0}^{i=I+1}$ be a set of uniformly distributed grid
points in $[-c_1,0]$ with $x_3^{i+1}-x_3^i=h$. Let $u_{l,i}^{(m,n)}$ be the
finite difference solution of
$u_l^{(m,n)}(x_3), l=1,2,3$ at the point $x_3=x_3^i$. The discrete finite
difference systems \eqref{odeu1} can be written as
\begin{equation}\label{dfu1}
\frac{u_{l,i-1}^{(m,n)}-2u_{l,i}^{(m,n)}+u_{l,i+1}^{(m,n)}}{h^2}+\Big{(}\kappa_1^2-\big{(}\frac{m\pi}{a}\big{)}^2-\big{(}\frac{n\pi}{b}\big{)}^2\Big{)}u_{l,i}^{(m,n)}=0,~ i=1,2,\cdots,I,
\end{equation}
where $(m,n)\in\mathbb{N}^2_l$, $l=1,2,3$.

Next, we consider the continuity conditions on $\Gamma_1=\{\boldsymbol
x\in\mathbb R^3|(x_1,x_2)\in[0,a]\times[0,b], x_3=-c_1\}$. By Maxwell's
equations, the tangential traces of the electromagnetic fields are continuous,
i.e.,
\[
 \nu\times \boldsymbol{E}_1=\nu\times\boldsymbol{E}_2,\quad  \nu\times
\boldsymbol{H}_1=\nu\times\boldsymbol{H}_2,
\]
the normal components of the electric and magnetic flux density are
continuous, i.e.,
\[
  \nu\cdot(\varepsilon_1 \boldsymbol{E}_1)=\nu \cdot (\varepsilon_2
\boldsymbol{E}_2),\quad \nu\cdot(\mu_0 \boldsymbol{H}_1)=\nu \cdot (\mu_0
\boldsymbol{H}_2).
\]
In addition, the electric field is divergence free, i.e.,
\[
 \nabla\cdot\boldsymbol{E}_1=\nabla\cdot\boldsymbol{E}_2=0.
\]
Componentwisely, the above continuity and divergence free conditions are
\begin{align}
\label{con11}&u_1(x_1,x_2,-c_1)=v_1(x_1,x_2,-c_1),\\
\label{con21}&u_2(x_1,x_2,-c_1)=v_2(x_1,x_2,-c_1),\\
\label{con31}&\varepsilon_1u_3(x_1,x_2,-c_1)=\varepsilon_2v_3(x_1,x_2,-c_1),\\
\label{con41}&\partial_{x_1}u_3(x_1,x_2,-c_1)-\partial_{x_3}u_1(x_1,x_2,-c_1)=\partial_{x_1}v_3(x_1,x_2,-c_1)-\partial_{x_3}v_1(x_1,x_2,-c_1),\\
\label{con51}&\partial_{x_2}u_3(x_1,x_2,-c_1)-\partial_{x_3}u_2(x_1,x_2,-c_1)=\partial_{x_2}v_3(x_1,x_2,-c_1)-\partial_{x_3}v_2(x_1,x_2,-c_1),\\
\label{con61}&\partial_{x_3}u_3(x_1,x_2,-c_1)=\partial_{x_3}v_3(x_1,x_2,-c_1).
\end{align}
Substituting \eqref{E1F}--\eqref{E2F} into \eqref{con11} and matching the
modes for the Fourier series expansions,  we obtain
\begin{equation*}
  u_1^{(m,n)}(-c_1)= v_1^{(m,n)}(-c_1),\quad (m,n)\in \mathbb{N}_1^2,
\end{equation*}
which implies
\begin{equation}\label{con12}
u_{1,0}^{(m,n)}=v_{1,J+1}^{(m,n)},\quad (m,n)\in \mathbb{N}_1^2.
\end{equation}

Similarly, we have from \eqref{con21}--\eqref{con31} that
\begin{equation}\label{con22}
u_{2,0}^{(m,n)}=v_{2,J+1}^{(m,n)},\quad (m,n)\in \mathbb{N}_2^2,
\end{equation}
and
\begin{equation}\label{con32}
\varepsilon_1 u_{3,0}^{(m,n)}=\varepsilon_2 v_{3,J+1}^{(m,n)}\quad (m,n)\in \mathbb{N}_3^2.
\end{equation}
Substituting \eqref{E1F}--\eqref{E2F} into \eqref{con41}, multiplying
the resulting equation by ${\rm cos}\big{(}\frac{m\pi
x_1}{a}\big{)}{\rm sin}\big{(}\frac{n\pi x_2}{b}\big{)}, (m,n)\in
\mathbb{N}^2_1$, and integrating over $\Gamma_1$, we obtain from the
orthogonality of the trigonometric functions that
\begin{equation*}
  \frac{\partial u_1^{(0,n)}(-c_1)}{\partial x_3}=\frac{\partial v_1^{(0,n)}(-c_1)}{\partial x_3},\quad n=1,2,\cdots,N,
\end{equation*}
and
\begin{equation*}
  u_3^{(m,n)}(-c_1)\frac{m\pi}{a}-\frac{\partial u_1^{(m,n)}(-c_1)}{\partial x_3}=v_3^{(m,n)}(-c_1)\frac{m\pi}{a}-\frac{\partial v_1^{(m,n)}(-c_1)}{\partial x_3}, \quad(m,n)\in\mathbb{N}^2_3.
\end{equation*}
Using the backward and forward finite difference schemes, we obtain
\begin{equation}\label{con42-1}
\frac{u_{1,1}^{(0,n)}-u_{1,0}^{(0,n)}}{h}=\frac{v_{1,J+1}^{(0,n)}-v_{1,J}^{(0,n)}}{h},\quad n=1,2,\cdots,N,
\end{equation}
and
\begin{equation}\label{con42-2}
\left(\frac{m\pi}{a}\right)u_{3,0}^{(m,n)}-\frac{u_{1,1}^{(m,n)}-u_{1,0}^{(m,n)}
}{h}=\left(\frac
{m\pi}{a}\right)v_{3,J+1}^{(m ,n)}
-\frac{v_{1,J+1}^{(m,n)}-v_{1,J}^{(m,n)}}{h},\quad(m,n)\in\mathbb{N} ^2_3.
\end{equation}
Combining \eqref{con42-1}--\eqref{con42-2}, \eqref{v1JJ1}, \eqref{con12} and
\eqref{con32} gives
\begin{equation}\label{con43-1}
\big{(}-1/r_{1,(JJ)}^{(0,n)}-2\big{)}u_{1,0}^{(0,n)}+u_{1,1}^{(0,n)}=0,\quad n=1,2,\cdots,N,
\end{equation}
and
\begin{equation}\label{con43-2}
\big{(}-1/r_{1,(JJ)}^{(m,n)}-2\big{)}u_{1,0}^{(m,n)}+\Big{(}\frac{\varepsilon_1}
{ \varepsilon_2}-1\Big{)}\left(\frac{m\pi
h}{a}\right)u_{3,0}^{(m,n)}+u_{1,1}^{(m,n)}=0,\quad(m,n)\in\mathbb{N}^2_3.
\end{equation}
For simplicity, let $u_{3,0}^{(0,n)}=u_{3,I+1}^{(0,n)}=0,~n=1,2,\cdots,N$ and
$u_{3,0}^{(m,0)}=u_{3,I+1}^{(m,0)}=0,~m=1,2,\cdots,M$ in the rest of this
section. Thus, \eqref{con43-1}--\eqref{con43-2} can be written uniformly as
\begin{equation}\label{con43}
\big{(}-1/r_{1,(JJ)}^{(m,n)}-2\big{)}u_{1,0}^{(m,n)}+\Big{(}\frac{\varepsilon_1}
{\varepsilon_2}-1\Big{)}\left(\frac{m\pi
h}{a}\right)u_{3,0}^{(m,n)}+u_{1,1}^{(m,n)}=0,\quad(m,n)\in\mathbb{N}^2_1.
\end{equation}
Similarly, based on the condition \eqref{con51}--\eqref{con61}, we get
\begin{equation}\label{con53}
\big{(}-1/r_{1,(JJ)}^{(m,n)}-2\big{)}u_{2,0}^{(m,n)}+\Big{(}\frac{\varepsilon_1}
{\varepsilon_2}-1\Big{)}\left(\frac{n\pi
h}{b}\right)u_{3,0}^{(m,n)}+u_{2,1}^{(m,n)}=0, \quad(m,n)\in\mathbb{N}^2_2,
\end{equation}
and
\begin{equation}\label{con63}
\Big{(}\big{(}-1/r_{2,(JJ)}^{(m,n)}-1\big{)}\frac{\varepsilon_1}{\varepsilon_2}-1\Big{)}u_{3,0}^{(m,n)}+u_{3,1}^{(m,n)}=0, \quad(m,n)\in\mathbb{N}^2_3.
\end{equation}

Define
\begin{equation*}
\boldsymbol{u}_l^{(m,n)}=\Big{(}u_{l,1}^{(m,n)},u_{l,2}^{(m,n)},\cdots,u_{l,I}^{
(m,n)}\Big{)}^\top,\quad l=1,2,3.
\end{equation*}
Using \eqref{con63}, we can rewrite the discrete system \eqref{dfu1} with $l=3$
in the matrix form
\begin{equation}\label{dfu3M}
\big{(}\boldsymbol{A}_4^{(m,n)}+\boldsymbol{D}_1^{(m,n)}\big{)}\boldsymbol{u}_3^{(m,n)}+\boldsymbol{a}_{I}u_{3,I+1}^{(m,n)}=0,~(m,n)\in\mathbb{N}^2_3,
\end{equation}
where
\begin{equation*}
  \boldsymbol{A}_4^{(m,n)}=\begin{pmatrix}
  1/\Big{(}(1/r_{2,(JJ)}^{(m,n)}+1)\frac{\varepsilon_1}{\varepsilon_2}+1)\Big{)}-2&1&&\\
  1&-2&1&\\
  &\ddots&\ddots&\ddots\\
  &&1&-2
  \end{pmatrix},\quad
  \boldsymbol{a}_{I}=\begin{pmatrix}
  0\\
  \vdots\\
   0\\
  1
  \end{pmatrix}.
\end{equation*}
and
\begin{equation*}
  \boldsymbol{D}_1^{(m,n)}=h^2\Big{(}\kappa_1^2-\big{(}\frac{m\pi}{a}\big{)}^2-\big{(}\frac{n\pi}{b}\big{)}^2\Big{)}\boldsymbol{I}_I.
\end{equation*}

Let
\begin{equation}\label{LUd13}
\boldsymbol{A}_4^{(m,n)}+\boldsymbol{D}_1^{(m,n)}=\boldsymbol{L}_4^{(m,n)}
\boldsymbol{U}_4^{(m,n)},\quad (m,n)\in\mathbb{N}^2_3
\end{equation}
be the LU-decomposition. It follows from
\eqref{dfu3M}--\eqref{LUd13} that
\begin{equation}\label{LUu31}
\boldsymbol{U}_4^{(m,n)}\boldsymbol{u}_3^{(m,n)}+\big{(}\boldsymbol{L}_4^{(m,n)}\big{)}^{-1}\boldsymbol{a}_{I}u_{3,I+1}^{(m,n)}=0,\quad (m,n)\in\mathbb{N}^2_3,
\end{equation}
where $\boldsymbol{U}_4^{(m,n)}=\big{(}r_{4,(pq)}^{(m,n)}\big{)}$. It follows
from \eqref{LUu31} that
\begin{equation}\label{LUu32}
\boldsymbol{u}_3^{(m,n)}=-\big{(}\boldsymbol{U}_4^{(m,n)}\big{)}^{-1}\big{(}\boldsymbol{L}_4^{(m,n)}\big{)}^{-1}\boldsymbol{a}_{I}u_{3,I+1}^{(m,n)},\quad (m,n)\in\mathbb{N}^2_3,
\end{equation}
It is clear to note that the first equation of the system \eqref{LUu32} is
\begin{equation}\label{u31T1}
u_{3,1}^{(m,n)}=(-1)^{2+I}\big({\rm
det}(\boldsymbol{U}_4^{(m,n)})\big)^{-1}u_{3,I+1}^{(m,n)},\quad
(m,n)\in\mathbb{N}^2_3.
\end{equation}

Let\begin{equation*}
  \boldsymbol{A}_3^{(m,n)}=\begin{pmatrix}
  1/(1/r_{1,(JJ)}^{(m,n)}+2)-2&1&&\\
  1&-2&1&\\
  &\ddots&\ddots&\ddots\\
  &&1&-2
  \end{pmatrix}_{I\times I},\quad
  \boldsymbol{a}_{1}^{(m,n)}=\begin{pmatrix}
  -1\\
  0\\
  \vdots\\
  0,
  \end{pmatrix}_{I\times 1},
\end{equation*}
where $(m,n)\in\mathbb{N}^2_1~{\rm or}~(m,n)\in\mathbb{N}^2_2$. Using
\eqref{con43}, \eqref{con63} and \eqref{u31T1}, we can write \eqref{dfu1} with
$l=1$ in the following matrix form:
\begin{equation}\label{dfu1M}
\big{(}\boldsymbol{A}_3^{(m,n)}+\boldsymbol{D}_1^{(m,n)}\big{)}\boldsymbol{u}_1^{(m,n)}+\boldsymbol{a}_{I}^{(m,n)}u_{1,I+1}^{(m,n)}=\boldsymbol{a}_{1}d_1^{(m,n)}u_{3,I+1}^{(m,n)},\quad(m,n)\in\mathbb{N}^2_1,
\end{equation}
where
\begin{equation*}
  d_1^{(0,n)}=0,\quad n=1,2,\cdots,N,
\end{equation*}
and
\begin{equation*}
  d_1^{(m,n)}=\left(\frac{(\frac{\varepsilon_1}{\varepsilon_2}-1)\frac{m\pi
h}{a}}{1/r_{1,(JJ)}^{(m,n)}+2}\right)\left(\frac{1}{(1/r_{2,(JJ)}^{(m,n)}
+1)\frac{\varepsilon_1}{\varepsilon_2}+1}\right)\left(\frac{(-1)^{2+I}}{{\rm
det}(\boldsymbol{U}_4^{(m ,n)})}\right), \quad(m,n)\in\mathbb{N}^2_3.
\end{equation*}

Similarly, we can rewrite \eqref{dfu1} with $l=2$ in the following matrix form:
\begin{equation}\label{dfu2M}
\big{(}\boldsymbol{A}_3^{(m,n)}+\boldsymbol{D}_1^{(m,n)}\big{)}\boldsymbol{u}_2^{(m,n)}+\boldsymbol{a}_{I}^{(m,n)}u_{2,I+1}^{(m,n)}=\boldsymbol{a}_{1}d_2^{(m,n)}u_{3,I+1}^{(m,n)},\quad(m,n)\in\mathbb{N}^2_2,
\end{equation}
where
\begin{equation*}
  d_2^{(m,0)}=0,\quad m=1,2,\cdots,M,
\end{equation*}
and
\begin{equation*}
  d_2^{(m,n)}=\left(\frac{(\frac{\varepsilon_1}{\varepsilon_2}-1)\frac{n\pi
h}{b}}{1/r_{1,(JJ)}^{(m,n)}+2}\right)\left(\frac{1}{(1/r_{2,(JJ)}^{(m,n)}
+1)\frac{\varepsilon_1}{\varepsilon_2}+1}\right)\left(\frac{(-1)^{2+I}}{
{\rm det}(\boldsymbol{U}_4^{(m ,n)})}\right), \quad(m,n)\in\mathbb{N}^2_3.
\end{equation*}
Let
\begin{equation}\label{LUd11}
  \boldsymbol{A}_3^{(m,n)}+\boldsymbol{D}_1^{(m,n)}=\boldsymbol{L}_3^{(m,n)}\boldsymbol{U}_3^{(m,n)},\quad (m,n)\in\mathbb{N}^2_1~{\rm or} ~(m,n)\in\mathbb{N}^2_2,
\end{equation}
be the LU-decomposition. It follows from \eqref{dfu1M}--\eqref{LUd11} that
\begin{equation}\label{LUu11}
\boldsymbol{U}_3^{(m,n)}\boldsymbol{u}_l^{(m,n)}+\big{(}\boldsymbol{L}_3^{(m,n)}\big{)}^{-1}\boldsymbol{a}_{I}u_{l,I+1}^{(m,n)}=\big{(}\boldsymbol{L}_3^{(m,n)}\big{)}^{-1}\boldsymbol{a}_{1}d_l^{(m,n)}u_{3,I+1}^{(m,n)},\quad (m,n)\in\mathbb{N}^2_l, l= 1,2,
\end{equation}
where $\boldsymbol{U}_3^{(m,n)}=\big{(}r_{3,(pq)}^{(m,n)}\big{)}$.

Combining the last equations of the systems \eqref{LUu11} and \eqref{LUu31}
gives
\begin{align*}
&r_{3,(II)}^{(m,n)}u_{l,I}+u_{l,I+1}=s_l^{(m,n)}u_{3,I+1}^{(m,n)},
\quad (m,n)\in\mathbb{N}^2_l,\mbox{ for } l=1,2, \\
&r_{4,(II)}^{(m,n)}u_{3,I}+u_{3,I+1}=0,\quad (m,n)\in\mathbb{N}^2_3,
\end{align*}
where $s_l^{(m,n)}=-\tilde{l}_{I1}^{(m,n)}d_l^{(m,n)}$, $l=1,2$, and
$\tilde{l}_{I1}^{(m,n)}$ is the $(I,1)$-th entry of
$\big{(}\boldsymbol{L}_3^{(m,n)}\big{)}^{-1}$.
We can write the above system in the matrix form
\begin{equation}\label{uTT1}
  \begin{pmatrix}
  \boldsymbol{R}_4&&\\
  &\boldsymbol{R}_5&\\
  &&\boldsymbol{R}_6
  \end{pmatrix}
    \begin{pmatrix}
  \boldsymbol{u}_{1,I}\\
  \boldsymbol{u}_{2,I}\\
  \boldsymbol{u}_{3,I}
  \end{pmatrix}
  +
  \begin{pmatrix}
  \hat{\boldsymbol{I}}_4&&-\mathcal{D}_1\\
  &\hat{\boldsymbol{I}}_5&-\mathcal{D}_2\\
  &&\hat{\boldsymbol{I}}_6
  \end{pmatrix}
    \begin{pmatrix}
  \boldsymbol{u}_{1,I+1}\\
  \boldsymbol{u}_{2,I+1}\\
  \boldsymbol{u}_{3,I+1}
  \end{pmatrix}
  =
 0,
\end{equation}
where
\begin{align*}
  \boldsymbol{R}_4 &={\rm diag}\big{(}r_{3,(II)}^{(m,n)}\big{)},
\quad(m,n)\in\mathbb{N}^2_1,\\
  \boldsymbol{R}_5 &={\rm diag}\big{(}r_{3,(II)}^{(m,n)}\big{)},
\quad(m,n)\in\mathbb{N}^2_2,\\
  \boldsymbol{R}_6 &={\rm diag}\big{(}r_{4,(II)}^{(m,n)}\big{)},
\quad(m,n)\in\mathbb{N}^2_3,
\end{align*}
the matrix $\mathcal{D}_1$ is the diagonal matrix ${\rm
diag}\big{(}s_1^{(m,n)}\big{)},~(m,n)\in\mathbb{N}^2_1$ by deleting the column
with respect to $m=0$, and the matrix $\mathcal{D}_2$ is the diagonal matrix
${\rm diag}\big{(}s_2^{(m,n)}\big{)},~(m,n)\in\mathbb{N}^2_2$ by deleting the
column with respect to $n=0$.\par
Similar to the homogeneous medium case, the TBC \eqref{ENU1} can be discretized as
\begin{equation}\label{uNUGa}
  \begin{pmatrix}
  \hat{\boldsymbol{I}}_1&&\\
  &\hat{\boldsymbol{I}}_2&\\
  &&\hat{\boldsymbol{I}}_3
  \end{pmatrix}
    \begin{pmatrix}
  \boldsymbol{u}_{1,I}\\
  \boldsymbol{u}_{2,I}\\
  \boldsymbol{u}_{3,I}
  \end{pmatrix}
  +
  \begin{pmatrix}
  -\hat{\boldsymbol{I}}_1-\boldsymbol{F}_1-\boldsymbol{H}_1&-\boldsymbol{G}_1&\boldsymbol{I}_1\\
  -\boldsymbol{H}_2&-\hat{\boldsymbol{I}}_2-\boldsymbol{F}_2-\boldsymbol{G}_2&\boldsymbol{I}_2\\
  \boldsymbol{F}_3&\boldsymbol{G}_3&-\hat{\boldsymbol{I}}_3
  \end{pmatrix}
    \begin{pmatrix}
  \boldsymbol{u}_{1,I+1}\\
  \boldsymbol{u}_{2,I+1}\\
  \boldsymbol{u}_{3,I+1}
  \end{pmatrix}
  =
  \begin{pmatrix}
  \boldsymbol{g}_1\\
  \boldsymbol{g}_2\\
  0
  \end{pmatrix}.
\end{equation}
Using \eqref{uTT1} and \eqref{uNUGa}, we obtain
\begin{equation}\label{uT1}
  \begin{pmatrix}
  -\hat{\boldsymbol{I}}_1-\boldsymbol{F}_1-\boldsymbol{H}_1-\boldsymbol{R}_4^{-1}&-\boldsymbol{G}_1&\boldsymbol{I}_1+\boldsymbol{R}_4^{-1}\mathcal{D}_1\\
  -\boldsymbol{H}_2&-\hat{\boldsymbol{I}}_2-\boldsymbol{F}_2-\boldsymbol{G}_2-\boldsymbol{R}_5^{-1}&\boldsymbol{I}_2+\boldsymbol{R}_5^{-1}\mathcal{D}_2\\
  \boldsymbol{F}_3&\boldsymbol{G}_3&-\hat{\boldsymbol{I}}_3-\boldsymbol{R}_6^{-1}
  \end{pmatrix}
    \begin{pmatrix}
  \boldsymbol{u}_{1,I+1}\\
  \boldsymbol{u}_{2,I+1}\\
  \boldsymbol{u}_{3,I+1}
  \end{pmatrix}
  =
  \begin{pmatrix}
  \boldsymbol{g}_1\\
  \boldsymbol{g}_2\\
  0
  \end{pmatrix}.
\end{equation}
The solution $\boldsymbol E_1=(u_1,u_2,u_3)$ on the open aperture $\Gamma$ can
be obtained by solving the linear system \eqref{uT1}.

\begin{remark}
For the cavity filled with a multi-layered medium (more than two layers),
a similar discretization can be developed for each layer, and similar discrete
continuity conditions can be deduced on the interface between every two
neighboring layers. As a result, a linear system similar to \eqref{uTT1} can be
obtained for the electric field in the first layer below the ground plane.
Consequently, we can get a linear system on the open aperture of the cavity by
using the linear system similar to \eqref{uTT1}--\eqref{uNUGa}.  The
solution $\boldsymbol{E}$ on the open aperture $\Gamma$ can be obtained by
solving the resulting system.
\end{remark}

\section{Evaluating singular integrals based on the FFT}\label{sec:FFT}

One of the key issues in the algorithm is how to evaluate efficiently and
accurately the singular integrals in \eqref{sig1}--\eqref{sig3}. Due to the lack
of closed form and the existence of singularity, direct numerical integration is
notoriously expensive. In this section, we propose an efficient algorithm to
evaluate these integrals based on the Fast Fourier Transform (FFT).
Specifically, we consider the evaluation of integrals
$\tilde{F}^{(m,n)}_{j,(k)}$, $\tilde{G}^{(m,n)}_{j,(k)}$, and $
\tilde{H}^{(m,n)}_{j,(k)}$, $j=1,2$ for $(m,n)\in\mathbb{N}^2, k\in\mathbb{N}^2,
\mathbb{N}^2=\{0,1,2,\cdots,M\}\times\{0,1,2,\cdots,N\}$. We refer to Appendix A
for the definition of $\tilde{F}^{(m,n)}_{2,(k)}$, $\tilde{G}^{(m,n)}_{2,(k)}$,
and $\tilde{H}^{(m,n)}_{2,(k)}$.

\subsection{Reduction of singularity}

It is easy to see that $\tilde{F}^{(m,n)}_{j,(k)}$, $j=1,2$ are weakly singular
integrals, while $\tilde{G}^{(m,n)}_{j,(k)}$, $\tilde{H}^{(m,n)}_{j,(k)}$,
$j=1,2$ include Cauchy type singular integrals. To make the computation easier,
we first apply the integration by parts to reduce the order of singularity
\begin{align*}
\tilde{G}^{(m,n)}_{1,(k)} &= \int_{\Gamma}\cos\left(\frac{m\pi
x_1}{a}\right)\sin\left(\frac{n\pi x_2}{b}\right)\partial_{x_2}
\left(\int_{\Gamma} \cos\left(\frac{k_1\pi y_1}{a}\right)\cos\left(\frac{k_2\pi
y_2}{b}\right)g(\boldsymbol{x},\boldsymbol{y}) ds_{\boldsymbol{y}} \right )
ds_{\boldsymbol{x}} \notag \\
&= \int_0^a \cos\left(\frac{m\pi x_1}{a}\right)\sin\left(\frac{n\pi
x_2}{b}\right) \left(\int_{\Gamma} \cos\left(\frac{k_1\pi
y_1}{a}\right)\cos\left(\frac{k_2\pi
y_2}{b}\right)g(\boldsymbol{x},\boldsymbol{y})  dy \right
)\Big|_{x_2=0}^{x_2=b}\ dx_1 \notag \\
&\quad -\int_{\Gamma}\partial_{x_2}\left(\cos\left(\frac{m\pi
x_1}{a}\right)\sin\left(\frac{n\pi x_2}{b}\right)\right) \left(\int_{\Gamma}
\cos\left(\frac{k_1\pi y_1}{a}\right)\cos\left(\frac{k_2\pi
y_2}{b}\right)g(\boldsymbol{x},\boldsymbol{y}) ds_{\boldsymbol{y}} \right )
ds_{\boldsymbol{x}} \notag \\
&= -\frac{n\pi}{b} \int_{\Gamma}\cos\left(\frac{m\pi
x_1}{a}\right)\cos\left(\frac{n\pi x_2}{b}\right) \left(\int_{\Gamma}
\cos\left(\frac{k_1\pi y_1}{a}\right)\cos\left(\frac{k_2\pi
y_2}{b}\right)g(\boldsymbol{x},\boldsymbol{y}) ds_{\boldsymbol{y}} \right )
ds_{\boldsymbol{x}}.
\end{align*}
Similar simplifications can be done for the integrals
$\tilde{G}^{(m,n)}_{2,(k)}$ and $\tilde{H}^{(m,n)}_{j,(k)}$, $j=1,2$. In the
end, we only need to consider evaluating the following three integrals:
\begin{align*}
I^{(m,n)}_{1,(k)}& = \int_{\Gamma}\cos\left(\frac{m\pi
x_1}{a}\right)\sin\left(\frac{n\pi x_2}{b}\right)\left(\int_{\Gamma}
\cos\left(\frac{k_1\pi y_1}{a}\right)\sin\left(\frac{k_2\pi
y_2}{b}\right)g(\boldsymbol{x},\boldsymbol{y}) ds_{\boldsymbol{y}} \right )
ds_{\boldsymbol{x}},\\
I^{(m,n)}_{2,(k)}& = \int_{\Gamma}\sin\left(\frac{m\pi
x_1}{a}\right)\cos\left(\frac{n\pi x_2}{b}\right)\left(\int_{\Gamma}
\sin\left(\frac{k_1\pi y_1}{a}\right)\cos\left(\frac{k_2\pi
y_2}{b}\right)g(\boldsymbol{x},\boldsymbol{y})ds_{\boldsymbol{y}} \right )
ds_{\boldsymbol{x}}, \\
I^{(m,n)}_{3,(k)}& = \int_{\Gamma}\cos\left(\frac{m\pi
x_1}{a}\right)\cos\left(\frac{n\pi x_2}{b}\right)\left(\int_{\Gamma}
\cos\left(\frac{k_1\pi y_1}{a}\right)\cos\left(\frac{k_2\pi
y_2}{b}\right)g(\boldsymbol{x},\boldsymbol{y}) ds_{\boldsymbol{y}} \right )
ds_{\boldsymbol{x}}.
\end{align*}
They belong to the same type of integrals, i.e.,
\begin{eqnarray*}
I^{(m,n)}_{(k)} =  \int_{\Gamma}\exp\left(\frac{m\pi x_1}{a}i\right)\exp\left(\frac{n\pi x_2}{b}i\right)\left( \int_{\Gamma} \exp\left(\frac{k_1\pi y_1}{a}i\right)\exp\left(\frac{k_2\pi y_2}{b}i\right)g(\boldsymbol{x},\boldsymbol{y}) ds_{\boldsymbol{y}}\right) ds_{\boldsymbol{x}}.
\end{eqnarray*}

Next, we propose a fast algorithm to evaluate $I^{(m,n)}_{(k)}$ by using the
FFT.

\subsection{Algorithm for $I^{(m,n)}_{(k)}$ based on the FFT}

Without loss of generality, we may assume $a \geq b$. To evaluate the integral $I^{(m,n)}_{(k)}$, we first consider evaluating the inner integral
\begin{eqnarray*}
I_{k}(x_1,x_2) = \int_{\Gamma} \exp\left(\frac{k_1\pi y_1}{a}i\right)\exp\left(\frac{k_2\pi y_2}{b}i\right)g(\boldsymbol{x},\boldsymbol{y}) dy
\end{eqnarray*}
for fixed  $k_1$, $k_2\in \mathbb{N}$ and $\boldsymbol{x}\in \Gamma$.

Define two functions:
\begin{eqnarray*}
f_{\mbox{Rect}}(x_1,x_2) = \begin{cases}
1, \quad \mbox{if } (x_1,x_2)\in \Gamma,\\
0, \quad \mbox{otherwise},
\end{cases}
\end{eqnarray*}
and
\begin{eqnarray*}
f_{\mbox{Circ}}(r) = \begin{cases}
1, \quad \mbox{if } r\le \sqrt{2}a,\\
0, \quad \mbox{otherwise}.
\end{cases}
\end{eqnarray*}
Then
\begin{eqnarray*}
I_{k}(x_1,x_2) = \int_{\mathbb{R}^2} \exp\left(\frac{k_1\pi
y_1}{a}i\right)\exp\left(\frac{k_2\pi
y_2}{b}i\right)f_{\mbox{Rect}}(y_1,y_2)g(\boldsymbol{x},\boldsymbol{y})
f_{\mbox{Circ}}(|\boldsymbol{x}-\boldsymbol{y}|) d\boldsymbol{y},\quad
\boldsymbol x\in\Gamma.
\end{eqnarray*}

Define $$F(\boldsymbol{y})= \exp\left(\frac{k_1\pi
y_1}{a}i\right)\exp\left(\frac{k_2\pi y_2}{b}i\right)f_{\mbox{Rect}}(y_1,y_2)$$
and $$G(\boldsymbol{x}-\boldsymbol{y})= g(\boldsymbol{x},\boldsymbol{y})
f_{\mbox{Circ}}(|\boldsymbol{x}-\boldsymbol{y}|). $$ Then
\begin{eqnarray*}
I_{k}(x_1,x_2) = \int_{\mathbb{R}^2} F(\boldsymbol{y})G(\boldsymbol{x}-\boldsymbol{y}) ds_{\boldsymbol{y}},
\end{eqnarray*}
which is a convolution and can be efficiently evaluated by using the FFT. Denote
by $\mathcal{F}(\cdot)$ the Fourier transform and $\mathcal{F}^{-1}(\cdot)$ the
inverse Fourier transform. Clearly, we have from the Fourier transformation
that
\begin{eqnarray*}
I_{k}(x_1,x_2) = \mathcal{F}^{-1}\left(\mathcal{F}(F)\cdot\mathcal{F}(G)\right).
\end{eqnarray*}
For $(j_1,j_2)\in \mathbb{N}^2$, it is easy to see that
\begin{align*}
\mathcal{F}(F)(j_1,j_2) &= \int_{\mathbb{R}^2}e^{-2\pi i(j_1y_1+
j_2y_2)}F(\boldsymbol{y})d\boldsymbol{y} \notag \\
&= \int_{\Gamma}e^{-2\pi i(j_1y_1+ j_2y_2)}\exp\left(\frac{k_1\pi
y_1}{a}i\right)\exp\left(\frac{k_2\pi y_2}{b}i\right)d\boldsymbol{y} \notag \\
&= \frac{(e^{-2\pi j_1a+ k_1\pi}-1)}{(-2\pi j_1+ k_1\pi/a)i}\frac{(e^{-2\pi
j_2b+ k_2\pi}-1)}{(-2\pi j_2+ k_2\pi/b)i}.
\end{align*}
Denote by $B$ the disk centered at the origin with radius $\sqrt{2}a$. The
following integral formula is convenient to evaluate $\mathcal{F}(G)(j_1,j_2)$:
\begin{align*}
\mathcal{F}(G)(j_1,j_2) &= \int_{\mathbb{R}^2}e^{-2\pi i(j_1y_1+
j_2y_2)}G(\boldsymbol{y})ds_{\boldsymbol{y}} \notag \\
&=\int_{B} e^{-2\pi i(j_1y_1+ j_2y_2)}
g(\boldsymbol{0},\boldsymbol{y})ds_{\boldsymbol{y}} \notag \\
&=\frac{1}{4\pi}\int_0^{\sqrt{2}a}\int_0^{2\pi} e^{-2\pi i(j_1 r\cos\theta +
j_2r\sin\theta)} \frac{e^{i\kappa_0r}}{r}rd\theta dr \notag \\
&=\frac{1}{4\pi} \int_0^{\sqrt{2}a} J_0\left(2\pi
\sqrt{j^2_1+j^2_2}r\right)e^{i\kappa_0 r} dr,
\end{align*}
where $J_0(\cdot)$ is the Bessel function of order zero.

Let $R =\sqrt{2}a$ and $c=2\pi \sqrt{j^2_1+j^2_2}$. Since there is no closed
form for the integral
\begin{eqnarray*}
I = \int_0^{R} J_0(c r)e^{i\kappa_0r} dr
\end{eqnarray*}
with $R>0$ and $c>0$, we need an algorithm to evaluate $I$ numerically.

We may assume $R$ and the wavenumber $\kappa_0$ are both $\mathcal{O}(1)$. Since $c=2\pi \sqrt{j^2_1+j^2_2}$ can be very large for $(j_1,j_2)\in \mathbb{N}^2$, in order to evaluate $I$ accurately, we consider two cases:
\begin{enumerate}
	\item Case 1: $j_1$ and $j_2$ are small, say, $\max\{|j_1|,|j_2|\}\le
10$, so that $c$ is $\mathcal{O}(1)$. In this case, direct integration by using
a high order Gaussian quadrature would efficiently evaluate the integral $I$.
	\item Case 2: $j_1$ and $j_2$ are large, in which case $c$ is large and $J_0(cr)$ is highly oscillatory. We can make use of the asymptotic formula
	\begin{eqnarray*}
	J_0(z) = \sqrt{\frac{2}{\pi z}}\left(\cos(z-\pi/4)+\frac{\sin(z-\pi/4)}{8z}+\mathcal{O}(\frac{1}{z^2})\right),
	\end{eqnarray*}
	which is quite accurate for $z\gg 1$ if we drop the reminder.
	Another useful formula is
	\begin{eqnarray*}
	\int_0^{\infty} J_0(c r)e^{i\kappa_0r} dr = \frac{1}{c^2-\kappa_0^2}, \mbox{ for } c>\kappa_0.
	\end{eqnarray*}
	Therefore,
	\begin{align*}
	I &= \frac{1}{c^2-\kappa_0^2} - \int_R^{\infty} J_0(c r)e^{i\kappa_0r}
dr \notag \\
	&= \frac{1}{c^2-\kappa_0^2} - \frac{1}{c}\int_{cR}^{\infty}
J_0(z)e^{i\kappa_0z/c} dz \notag \\
	&\approx \frac{1}{c^2-\kappa_0^2} -
\frac{1}{c}\sqrt{\frac{2}{\pi}}\int_{cR}^{\infty}
\left(\frac{\cos(z-\pi/4)}{\sqrt{z}}+\frac{\sin(z-\pi/4)}{8z^{3/2}} \right)
e^{i\kappa_0z/c} dz \notag \\
	&= \frac{1}{c^2-\kappa_0^2} -
\frac{1}{c}\sqrt{\frac{2}{\pi}}\int_{cR}^{\infty}
\left(\frac{e^{(z-\pi/4)i}+e^{(\pi/4-z)i}}{2\sqrt{z}}+\frac{e^{(z-\pi/4)i}-e^{
(z-\pi/4)i}}{2i8z^{3/2}} \right) e^{i\kappa_0z/c} dz \notag.
	\end{align*}
	In other words, we have to evaluate these two kinds of integrals
	\begin{eqnarray*}
	\int_{R_0}^\infty \frac{e^{p z i}}{\sqrt{z}}dz, \quad \int_{R_0}^\infty \frac{e^{q z i}}{z^{3/2}}dz,
	\end{eqnarray*}
	where $p,q\in \mathbb{R}$ and $R_0\gg 1$. They belong to the same type
of integrals. In fact, we obtain from the integration by parts that
	\begin{eqnarray*}
	\int_{R_0}^\infty \frac{e^{q z i}}{z^{3/2}}dz = 2\frac{e^{q {R_0} i}}{\sqrt{R_0}}+2qi\int_{R_0}^\infty \frac{e^{q z i}}{\sqrt{z}}dz.
	\end{eqnarray*}
	On the other hand,
	\begin{eqnarray*}
	\int_{R_0}^\infty \frac{e^{p z i}}{\sqrt{z}}dz =\frac{\sqrt{\pi}}{2p}\left(1+i-2 \mbox{Fresnelc}(\sqrt{2pR_0/\pi})-2i \mbox{Fresnels}(\sqrt{2pR_0/\pi})\right),
	\end{eqnarray*}
	where $\mbox{Fresnelc}(\cdot)$ and  $\mbox{Fresnels}(\cdot)$ are Fresnel
cosine  and sine integrals, respectively. To efficiently evaluate them, we make
use of the following asymptotic expansions for $z\gg 1$:
	\begin{eqnarray*}
	\mbox{Fresnelc}(z) &=& \frac{1}{2}+f(z)\sin(\frac{1}{2}\pi z^2)-g(z) \cos(\frac{1}{2}\pi z^2) ,\\
	\mbox{Fresnels}(z) &=& \frac{1}{2}-f(z)\cos(\frac{1}{2}\pi z^2)-g(z) \sin(\frac{1}{2}\pi z^2),
	\end{eqnarray*}
	where
	\begin{eqnarray*}
	f(z) &=& \frac{1}{\pi z}\left(1-\frac{3}{(\pi z^2)^2}+\mathcal{O}(\frac{1}{z^8})\right), \\
	g(z) &=& \frac{1}{\pi^2 z^3}\left(1-\frac{15}{(\pi z^2)^2}+\mathcal{O}(\frac{1}{z^8})\right).
	\end{eqnarray*}
\end{enumerate}

Combining all the ingredients above, we are able to efficiently evaluate the
inner integral $I_k(x_1,x_2)$. Once $I_k(x_1,x_2)$ is available, for the outer
integral with respect to $\boldsymbol{x}$, we simply use the trapezoidal rule,
in which case the FFT can also be directly applied.

\section{Implementation and complexity}

Our algorithm is extremely efficient in terms of computational cost. A detailed
analysis on the  computational complexity of Algorithm I for the
electromagnetic scattering by an open rectangular cavity filled with a
homogeneous medium and Algorithm II for the electromagnetic scattering by
an open rectangular cavity filled with a layered medium.

\begin{table}[ht]
\begin{tabular}{cp{.8\textwidth}}
\toprule
\multicolumn{2}{l}{Algorithm I:\quad Electromagnetic
scattering by a homogeneous cavity.} \\
\midrule
Step 1 &  Generate the matrices
$\boldsymbol{F}_i,\boldsymbol{G}_i,\boldsymbol{H}_i,i=1,2$ and the vectors
$\boldsymbol{g}_i, i=1,2$; \\
Step 2 & Calculate the LU decomposition to get
$\boldsymbol{U}_i^{(m,n)},i=1,2$ and  $\boldsymbol{R}_i^{-1},i=1,2,3$ by using
the forward Gaussian elimination with a row partial pivoting; \\
Step 3 &  Solve the system \eqref{EJ1} for
$\boldsymbol{E}_{i,J+1},i=1,2,3$. \\
\bottomrule
\end{tabular}
\end{table}

\begin{table}[ht]
\begin{tabular}{cp{.8\textwidth}}
\toprule
\multicolumn{2}{l}{Algorithm II:\quad Electromagnetic
scattering by a layered cavity.} \\
\midrule
Step 1 &  Generate the matrices
$\boldsymbol{F}_i,\boldsymbol{G}_i,\boldsymbol{H}_i,i=1,2$ and the vectors
$\boldsymbol{g}_i, i=1,2$; \\
Step 2 &  Calculate the LU decomposition to get
$\boldsymbol{U}_i^{(m,n)},i=1,2$ and  $\boldsymbol{R}_i^{-1},i=1,2,3$ by using
the forward Gaussian elimination with a row partial pivoting. Further, calculate
the LU decomposition to get
$\boldsymbol{U}_i^{(m,n)},\boldsymbol{U}_i^{(m,n)},i=3,4$ $\mathcal{D}_i,i=1,2$
and  $\boldsymbol{R}_i^{-1},i=4,5,6$; \\
Step 3 &  Solve the system \eqref{uT1} for $\boldsymbol{E}$ on
the open aperture $\Gamma$.\\
\bottomrule
\end{tabular}
\end{table}

The cost for each step is presented in Table \ref{Tab1}. In Step 1, one
needs to calculate the singular integrals to generate the matrices
$\boldsymbol{F}_i,\boldsymbol{G}_i,\boldsymbol{H}_i,i=1,2$. As shown in Section
\ref{sec:FFT}, we evaluate the singular integrals based on FFT, which requires
only $MN(MN\log(MN)+MN\log(MN))$ complex operations for all the singular
integrals. Hence the overall cost of Step 1 is  $O(M^2N^2\log(MN))$. In Step
2 for Algorithm I, we need to calculate the LU decomposition for
$\boldsymbol{A}_1+\boldsymbol{D}^{(m,n)},~(m,n)\in\mathbb{N}^2_1\cup\mathbb{N}
^2_2$ and $\boldsymbol{A}_2+\boldsymbol{D}^{(m,n)},~(m,n)\in\mathbb{N}^2_3$. By
noting the tridiagonal structure of these matrices, only $3J(M+1)(N+1)+3JMN$
complex operations are needed.  In Step 2 for Algorithm II, the cost for
calculating the LU decomposition in the bottom layer is $3J(M+1)(N+1)+3JMN$, and
the cost for calculating the LU decomposition in the top layer is
$5I(M+1)(N+1)+3IMN$. In Algorithms I and II, we need
to solve the interface system \eqref{EJ1} and \eqref{uT1}, respectively. We
point out that a direct method, such as the Gaussian elimination scheme,
requires $(3MN+2N)^3/3$ complex operations, which is not efficient. In order to
solve the interface system effectively, we may need the effective iterative
solver. The efficiency of the iterative algorithm for the interface system
depends upon many factors, such as the complicated transparent boundary
condition, the regularity of solution, the eigenvalue distribution and the
condition numbers of the coefficient matrix. We will carry out the related work
in the follow-up work.

\begin{table}[ht]
\caption{The computational complexity of Algorithms I and II.}\label{Tab1}
\begin{tabular}{ccccc}
\hline
\hline
Step && Homogeneous cavity && Layered cavity\\
\hline
1 &&$O(M^2N^2\log(MN))$&&$O(M^2N^2\log(MN))$\\
2 &&$3J((M+1)(N+1)+MN)$&&$3J(M+1)(N+1)+3JMN$\\
&& &&$+5I(M+1)(N+1)+3IMN$\\
3 &&$(3MN+2N)^3/3$&&$(3MN+2N)^3/3$\\
\hline
\hline
\end{tabular}
\end{table}

\section{Numerical experiments}

In this section, several numerical examples are presented to demonstrate the the
performance of the proposed method. Throughout all the examples, the incident
wave
\begin{equation*}
\boldsymbol{E}^{\rm inc}(\boldsymbol{x})=(\cos \alpha
\hat{\boldsymbol{\theta}}+\sin\alpha \hat{\boldsymbol{\phi}})e^{i\kappa_0
\boldsymbol{d} r},
\end{equation*}
where $\alpha$ is the polarization angle, $\hat{\theta}$ and $\hat{\phi}$ are
the standard unit vectors in the spherical coordinates, and $\boldsymbol d$ is
the incident direction given by
\begin{equation*}
  \boldsymbol{d}=-(\sin \theta\cos \phi,\sin \theta\sin \phi, \cos \theta).
\end{equation*}
The wavenumber $\kappa_0=2\pi$. The incident angle $\phi=0$ so that we focus on
the $xz$-plane.

The physical quantity of interest associated with the cavity scattering is the
radar cross section (RCS), which measures the detectability of a target by a
radar system \cite{J02}. When the incident angle and the observation angle
are the same, the RCS is called the backscatter RCS. The specific formulas can
be found in \cite{J91} for the RCS of the three-dimensional cavity-backed
apertures.

Our fast algorithm is mainly validated and compared with the adaptive finite
element PML method. The fast algorithm is carried out by a laptop with
Intel(R) Core(TM) i5-2430M CPU @ 2.40GHz. The implementation of the adaptive
finite element PML method is based on parallel hierarchical grid (PHG)
\cite{phg,ZZC16}, which is a toolbox for developing parallel adaptive finite
element programs on unstructured tetrahedral meshes. The linear system resulted
from the finite element discretization is solved by MUMPS (MUltifrontal
Massively Parallel Sparse direct Solver) \cite{mumps}, which is a general
purpose library for the direct solution of large linear systems. The
computation is done on the high performance computers of State Key Laboratory
of Scientific and Engineering Computing, Chinese Academy of Sciences, in which
each node has 2 Intel Xeon Gold 6140 CPUs (2.3 GHz, 18 cores) and 192 GB memory
and a 100 GB EDR Infiniband network is used for data communication between
nodes. We solve the finite element problem for each $\theta$ with one node (36
cores). The maximum number of degrees of freedom (DoFs) on the mesh are between
2,000,000 and 3,000,000. The running time (CPU time/cores) is 5 to 10 minutes.
By choosing the increment of $\theta$ as $\Delta\theta=0.5^\circ$, the
finite element problem is solved 100 times in Example 1, and 180 times in
Examples 2 and 3.

When presenting the numerical results, we use the following notations:
\begin{itemize}
	\item $M$, $N$: Number of modes for the Fourier expansions in the $x_1$
and $x_2$ directions, respectively.
	\item $J$: Number of partition points along the $x_3$ direction. For
a two-layered medium, another variable $I$ is used.
	\item $T_{\rm singular}$: Amount of time in seconds required to evaluate
the singular integrals.
	\item $T_{\rm assemble}$: Amount of time in seconds required to assemble
the matrix.
	\item $T_{\rm solve}$: Amount of time in seconds required to solve the
linear system.
	\item $T_{\rm RCS}$: Amount of time in seconds required to calculate the
RCS.
\end{itemize}

\subsection{Example 1}

\begin{figure}
\centering
\includegraphics[width=0.45\textwidth]{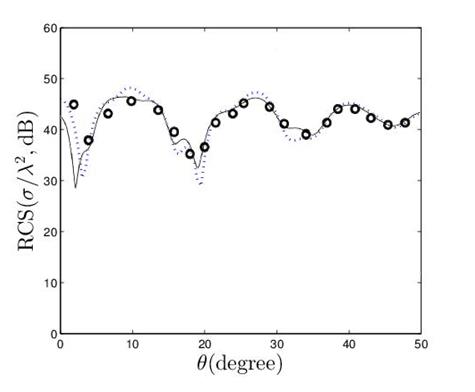}\quad
\includegraphics[width=0.45\textwidth]{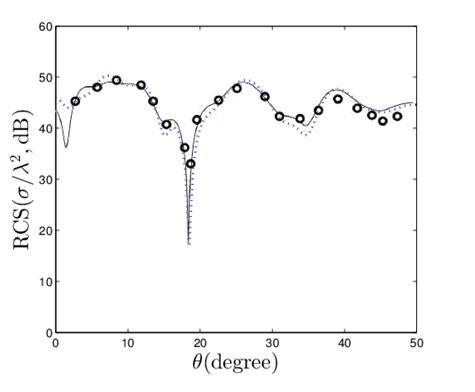}
\caption{Example 1: the backscatter RCS of the cavity with size $a=b=10\lambda$
and $c=30\lambda$. The dashed line is the RCS calculated by
our fast algorithm, the solid line is the RCS calculated by the mode matching
method presented in \cite{BGLZ12}, and the circle is the RCS calculated by the
modal approach presented in \cite{LLC89}.}
\label{fig:3}
\end{figure}

\begin{figure}
\centering
\includegraphics[width=0.45\textwidth]{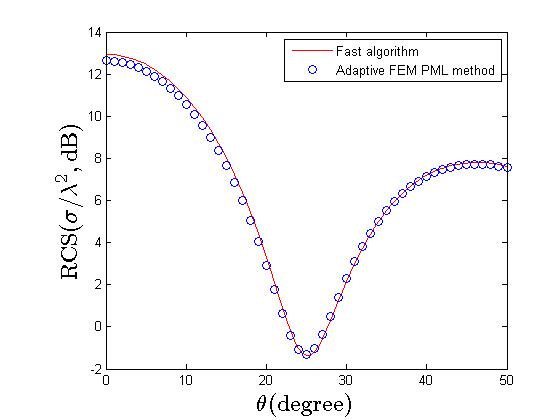}\quad
\includegraphics[width=0.45\textwidth]{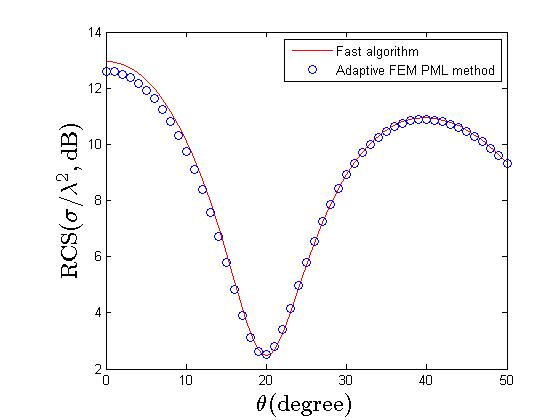}
\caption{Example 1: the backscatter RCS of the cavity with size $a=b=\lambda$
and $c=3\lambda$. }
\label{fig:4}
\end{figure}

In this example, we consider the cavity filled with a homogeneous medium.
First, the backscatter
RCS of the cavity with size $a=b=10\lambda$ and $c=30\lambda$ is calculated.
The RCS of $\hat{\theta}\hat{\theta}$ and $\hat{\phi}\hat{\phi}$ polarizations
are shown in Figure \ref{fig:3} for various incident angle $\theta$. The
numerical results show excellent agreement with the calculations by the mode
matching method presented in \cite{BGLZ12} and the modal approach presented in
\cite{LLC89}. The detailed computational time is given in Table \ref{Tab0}. Most
of the time is spent on the evaluation of singular integrals. However, we only
need to compute them once for different incident angles. In addition, most of
the applications only require a small number of modes to resolve the field.
Next, the backscatter radar cross section of a cavity with size $a=b=\lambda$
and $c=3\lambda$ is calculated by the fast algorithm and the adaptive finite
element PLM method. Figure \ref{fig:4} shows the RCS versus $\theta$ for
$\hat{\theta}\hat{\theta}$ and $\hat{\phi}\hat{\phi}$ polarizations. The
backscatter RCS is shown as red solid lines and blue circles for the fast
algorithm and adaptive PML method, respectively. It is clear to note that the
results obtained by both methods are consistent with each other. Detailed computational time is given in Table \ref{Tab2}.

\begin{table}[ht]
	\caption{Example 1: The time of the computation for the cavity size $a=b=\lambda,c=3\lambda$ with $\alpha=0,\theta=\pi/6$.}\label{Tab0}
	\begin{tabular}{cccccc}
		\hline
		\hline
		$M$, $N$ & $J$ & $T_{\rm singular}$ & $T_{\rm assemble}$ &
$T_{\rm solve}$ & $T_{\rm RCS}$ \\
		\hline
		$M=N=21$&1000&119.697858&0.307505&0.286226& 0.006255\\
		\hline
		$M=N=3$&1000&4.100736&0.002453&0.000072& 0.000715\\
		\hline
		$M=N=3$&600&4.013451& 0.002015&0.000066&0.000691\\
		\hline
		\hline
	\end{tabular}
\end{table}

\begin{table}[ht]
	\caption{Example 1: The time of the computation for the cavity size $a=b=10\lambda,c=30\lambda$ with $\alpha=0,\theta=\pi/6$.}\label{Tab2}
	\begin{tabular}{cccccc}
		\hline
		\hline
		$M$, $N$ & $J$ & $T_{\rm singular}$ & $T_{\rm assemble}$ &
$T_{\rm solve}$ & $T_{\rm RCS}$ \\
		\hline
		$M=N=21$&1000&212.170604&0.298627&0.335862&0.006811\\
		\hline
		$M=N=15$&1000&155.201708&0.062623&0.058558&0.003075\\
		\hline
		$M=N=15$&1500&155.360674&0.063114&0.038757& 0.006100\\
		\hline
		\hline
	\end{tabular}
\end{table}

\subsection{Example 2}

\begin{figure}
\centering
\includegraphics[width=0.45\textwidth]{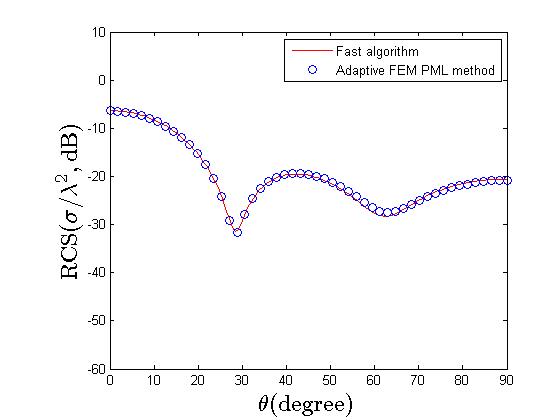}\quad
\includegraphics[width=0.45\textwidth]{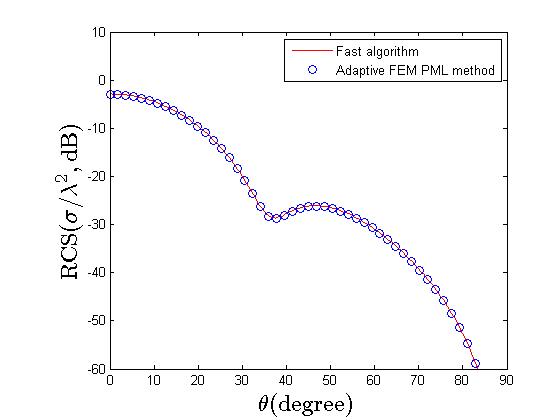}
\caption{Example 2: the backscatter RCS of the cavity with size $a=\lambda$ and
$b=c=0.25\lambda$. }
\label{fig:5}
\end{figure}

In this example, we consider the cavity filled with a material having a
relative permittivity $\epsilon_{\rm r}=7+1.5 {\,\rm i}$ and a constant magnetic
permeability $\mu=1$. The backscatter RCS of the cavity with size
$a=\lambda$,  and $b=c=0.25\lambda$ is calculated.  The RCS of
$\hat{\theta}\hat{\theta}$ and $\hat{\phi}\hat{\phi}$ polarizations are shown in
Figure \ref{fig:5} for various incident angle $\theta$. The results based on
the fast algorithm and the adaptive PML method are again in excellent agreement.
Detailed computational time is given in Table \ref{Tab3}. Again, the cost is
dominated by the evaluation of singular integrals.

\begin{table}[ht]
	\caption{Example 2: The time of the computation for the cavity size $a=\lambda,b=c=0.25\lambda$ with $\alpha=0,\theta=\pi/2$.}\label{Tab3}
	\begin{tabular}{cccccc}
		\hline
		\hline
		$M$, $N$ & $J$ & $T_{\rm singular}$ & $T_{\rm assemble}$ &
$T_{\rm solve}$ & $T_{\rm RCS}$ \\
		\hline
		$M=N=15$&1000&43.997846&0.065300&0.051353&0.003815\\
		\hline
		$M=N=3$&1000&3.037277&0.001508&0.000068&0.001106\\
		\hline
		$M=N=3$&100& 3.047187& 0.001897&0.000098&0.001130\\
		\hline
		\hline
	\end{tabular}
\end{table}

\subsection{Example 3}

\begin{figure}
\centering
\includegraphics[width=0.45\textwidth]{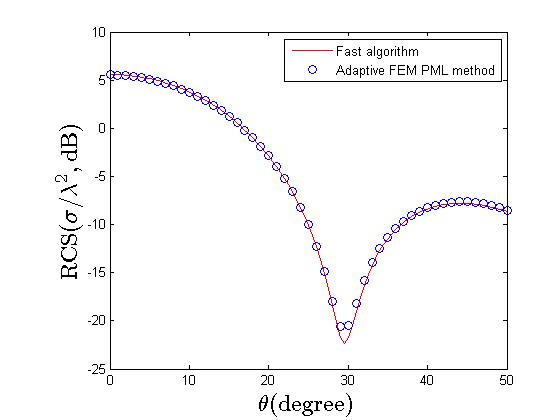}\quad
\includegraphics[width=0.45\textwidth]{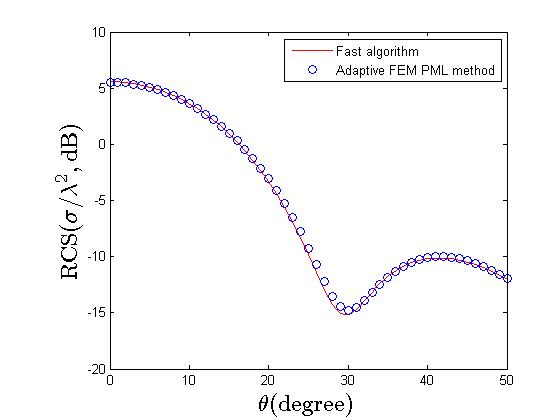}
\caption{Example 3: The backscatter RCS of the cavity with size $a=b=\lambda$,
$c_1=\lambda$ and $c_2=2\lambda$. }
\label{fig:6}
\end{figure}

This example is concerned with the cavity filled with a two-layer material.
The cavity size is $a=b=\lambda$ and $c=3\lambda$. The top and  bottom layer
materials have parameters $\epsilon_{\rm r}=7+1.5 {\,\rm i}$ and $\epsilon_{\rm
r}=3+0.05 {\,\rm i}$, respectively. The thickness of the top material and the
bottom material are $c_1=\lambda$ and $c_2=2\lambda$, respectively. The
backscatter RCS of $\hat{\theta}\hat{\theta}$ and $\hat{\phi}\hat{\phi}$
polarizations are shown in Figure \ref{fig:6} for various incident angle
$\theta$. Once again, both methods are consistent with each other very well.
Detailed computational time is given in Table \ref{Tab4}. It can be seen that
the total computational time is less than three minutes by using our fast
algorithm.

\begin{table}[ht]
	\caption{Example 3: The time of the computation for the cavity size $a=b=\lambda,c_1=\lambda,c_2=2\lambda$ with $\alpha=0,\theta=\pi/2$.}\label{Tab4}
	\begin{tabular}{cccccc}
		\hline
		\hline
		$M$, $N$ & $J$, $I$ & $T_{\rm singular}$ & $T_{\rm assemble}$ &
$T_{\rm solve}$ & $T_{\rm RCS}$ \\
		\hline
		$M=N=21$&$J=1000,I=2000$&130.499710&0.329729&0.341623&0.006732\\
		\hline
		$M=N=3$&$J=1000,I=2000$&4.075719&0.002566&0.000097& 0.001600\\
		\hline
		$M=N=3$&$J=100,I=200$& 4.194746&0.002448&0.000075&0.001176\\
		\hline
		\hline
	\end{tabular}
\end{table}

\section{Conclusion}

In this paper, we have presented a fast algorithm for the electromagnetic
scattering from three dimensional open rectangular cavities. Based on the
Fourier series expansions in the horizontal directions and the Gaussian
elimination along the vertical direction, the fast algorithm reduces the global
system to an interface system on the open aperture only. We also propose an
efficient algorithm to evaluate the singular integrals on the aperture based on
FFT. The whole algorithm enjoys the advantage of the low computational cost by
solving only the coefficients of the modes of the Fourier series. Moreover, our
fast algorithm has the capability of handling large cavities or high wave
numbers. This work provides a viable alternative to the current efforts of
designing sophisticated
basis functions for solving Maxwell equations with high wave numbers or large
cavities. A possible future work is to extend our fast algorithm to the optimal
design problems and inverse problems. It is of particular interest in designing the shape and
composition of a layered cavity to minimize the RCS~\cite{BJL1,BJL2}. Computationally, the design
problem can be challenging because of the need of solving
the scattering problem repeatedly. The fast algorithm presented here certainly would
provide an efficient and accurate numerical tool for these problems.

\section*{Appendix A. Notations}

In the appendix, we list the expressions of the entries for the vectors
used in \eqref{ENU14}.

The definitions of $\boldsymbol{F}_{1}$, $\boldsymbol{G}_{1}$ and
$\boldsymbol{H}_{1}$ are given by
	\begin{equation*}
	\boldsymbol{F}_{1}:=\begin{pmatrix}
	\boldsymbol{F}_{1}^{(0)}\\
	\boldsymbol{F}_{1}^{(1)}\\
	\vdots\\
	\boldsymbol{F}_{1}^{(M)}
	\end{pmatrix},\quad
	\boldsymbol{G}_{1}:=\begin{pmatrix}
	\boldsymbol{G}_{1}^{(0)}\\
	\boldsymbol{G}_{1}^{(1)}\\
	\vdots\\
	\boldsymbol{G}_{1}^{(M)}
	\end{pmatrix},\quad
	\boldsymbol{H}_{1}:=\begin{pmatrix}
	\boldsymbol{H}_{1}^{(0)}\\
	\boldsymbol{H}_{1}^{(1)}\\
	\vdots\\
	\boldsymbol{H}_{1}^{(M)}
	\end{pmatrix},
	\end{equation*}
with	
	\begin{equation*}
	\boldsymbol{F}_{1}^{(m)}:=\begin{pmatrix}
	\boldsymbol{F}_{1}^{(m,1)}\\
	\boldsymbol{F}_{1}^{(m,2)}\\
	\vdots\\
	\boldsymbol{F}_{1}^{(m,N)}
	\end{pmatrix},\quad
	\boldsymbol{G}_{1}^{(m)}:=\begin{pmatrix}
	\boldsymbol{G}_{1}^{(m,1)}\\
	\boldsymbol{G}_{1}^{(m,2)}\\
	\vdots\\
	\boldsymbol{G}_{1}^{(m,N)}
	\end{pmatrix},\quad
	\boldsymbol{H}_{1}^{(m)}:=\begin{pmatrix}
	\boldsymbol{H}_{1}^{(m,1)}\\
	\boldsymbol{H}_{1}^{(m,2)}\\
	\vdots\\
	\boldsymbol{H}_{1}^{(m,N)}
	\end{pmatrix},
	\end{equation*}
where
\begin{align*}
\boldsymbol{F}_{1}^{(m,n)}:=\begin{pmatrix}
\boldsymbol{F}_{1,(0)}^{(m,n)}&\boldsymbol{F}_{1,(1)}^{(m,n)}&\cdots&\boldsymbol{F}_{1,(M)}^{(m,n)}
\end{pmatrix},\quad
&\boldsymbol{F}_{1,(k_1)}^{(m,n)}:=\begin{pmatrix}
F_{1,(k_1,1)}^{(m,n)}&F_{1,(k_1,2)}^{(m,n)}&\cdots&F_{1,(k_1,N)}^{(m,n)}
\end{pmatrix},\\
\boldsymbol{G}_{1}^{(m,n)}:=\begin{pmatrix}
\boldsymbol{G}_{1,(1)}^{(m,n)}&\boldsymbol{G}_{1,(2)}^{(m,n)}&\cdots&\boldsymbol{G}_{1,(M)}^{(m,n)}
\end{pmatrix},\quad
&\boldsymbol{G}_{1,(k_1)}^{(m,n)}:=\begin{pmatrix}
G_{1,(k_1,0)}^{(m,n)}&G_{1,(k_1,1)}^{(m,n)}&\cdots&G_{1,(k_1,N)}^{(m,n)}
\end{pmatrix},\\
\boldsymbol{H}_{1}^{(m,n)}:=\begin{pmatrix}
\boldsymbol{H}_{1,(0)}^{(m,n)}&\boldsymbol{H}_{1,(1)}^{(m,n)}&\cdots&\boldsymbol{H}_{1,(M)}^{(m,n)}
\end{pmatrix},\quad
&\boldsymbol{H}_{1,(k_1)}^{(m,n)}:=\begin{pmatrix}
H_{1,(k_1,1)}^{(m,n)}&H_{1,(k_1,2)}^{(m,n)}&\cdots&H_{1,(k_1,N)}^{(m,n)}
\end{pmatrix}.
\end{align*}

The definitions of $\boldsymbol{g}_1$ and $\boldsymbol{E}_{l,j}$ with $l=1,2,3$,
$0\le j\le J+1$ are given by
\begin{align*}
&\boldsymbol{g}_1:=\begin{pmatrix}
\boldsymbol{g}_{1}^{(0)}&\boldsymbol{g}_{1}^{(1)}&
\hdots&\boldsymbol{g}_{1}^{(M)}
\end{pmatrix}^\top, \quad \boldsymbol{g}_{1}^{(m)}:=\begin{pmatrix}
g_{1}^{(m,1)}&g_{1}^{(m,2)}&\hdots&g_{1}^{(m,N)}
\end{pmatrix},\\
& \boldsymbol{E}_{1,j}:=\begin{pmatrix}
  \boldsymbol{E}_{1,j}^{(0)}&\boldsymbol{E}_{1,j}^{(1)}&\hdots&\boldsymbol{E}_{1,j}^{(M)}
  \end{pmatrix}^\top, \quad \boldsymbol{E}_{1,j}^{(m)}:=\begin{pmatrix}
  E_{1,j}^{(m,1)}&E_{1,j}^{(m,2)}&\hdots&E_{1,j}^{(m,N)}
  \end{pmatrix},\\
  &\boldsymbol{E}_{2,j}:=\begin{pmatrix}
  \boldsymbol{E}_{2,j}^{(1)}&\boldsymbol{E}_{2,j}^{(2)}&\hdots&\boldsymbol{E}_{2,j}^{(M)}
  \end{pmatrix}^\top, \quad
  \boldsymbol{E}_{2,j}^{(m)}:=\begin{pmatrix}
  E_{2,j}^{(m,0)}&E_{2,j}^{(m,1)}&\hdots&E_{2,j}^{(m,N)}
  \end{pmatrix},\\
  &\boldsymbol{E}_{3,j}:=\begin{pmatrix}
  \boldsymbol{E}_{3,j}^{(1)}&\boldsymbol{E}_{3,j}^{(2)}&
  \hdots&\boldsymbol{E}_{3,j}^{(M)}
  \end{pmatrix}^\top, \quad \boldsymbol{E}_{3,j}^{(m)}:=\begin{pmatrix}
  E_{3,j}^{(m,1)}&E_{3,j}^{(m,2)}&
  \hdots&E_{3,j}^{(m,N)}
  \end{pmatrix}.
\end{align*}

The definitions of $\boldsymbol{F}_{2}$, $\boldsymbol{G}_{2}$ and
$\boldsymbol{H}_{2}$ are given by
\begin{equation*}
\boldsymbol{F}_{2}:=\begin{pmatrix}
\boldsymbol{F}_{2}^{(1)}\\
\boldsymbol{F}_{2}^{(2)}\\
\vdots\\
\boldsymbol{F}_{2}^{(M)}
\end{pmatrix}, \quad
\boldsymbol{G}_{2}:=\begin{pmatrix}
\boldsymbol{G}_{2}^{(1)}\\
\boldsymbol{G}_{2}^{(2)}\\
\vdots\\
\boldsymbol{G}_{2}^{(M)}
\end{pmatrix}, \quad
\boldsymbol{H}_{2}:=\begin{pmatrix}
\boldsymbol{H}_{2}^{(1)}\\
\boldsymbol{H}_{2}^{(2)}\\
\vdots\\
\boldsymbol{H}_{2}^{(M)}
\end{pmatrix},
\end{equation*}
with
\begin{equation*}
\boldsymbol{F}_{2}^{(m)}:=\begin{pmatrix}
\boldsymbol{F}_{2}^{(m,0)}\\
\boldsymbol{F}_{2}^{(m,1)}\\
\vdots\\
\boldsymbol{F}_{2}^{(m,N)}
\end{pmatrix}, \quad
\boldsymbol{G}_{2}^{(m)}:=\begin{pmatrix}
\boldsymbol{G}_{2}^{(m,0)}\\
\boldsymbol{G}_{2}^{(m,1)}\\
\vdots\\
\boldsymbol{G}_{2}^{(m,N)}
\end{pmatrix}, \quad
\boldsymbol{H}_{2}^{(m)}:=\begin{pmatrix}
\boldsymbol{H}_{2}^{(m,0)}\\
\boldsymbol{H}_{2}^{(m,1)}\\
\vdots\\
\boldsymbol{H}_{2}^{(m,N)}
\end{pmatrix},
\end{equation*}
where
\begin{align*}
&\boldsymbol{F}_{2}^{(m,n)}:=\begin{pmatrix}
\boldsymbol{F}_{2,(1)}^{(m,n)}&\boldsymbol{F}_{2,(2)}^{(m,n)}&\cdots&\boldsymbol{F}_{2,(M)}^{(m,n)}
\end{pmatrix}, \quad
\boldsymbol{F}_{2,(k_1)}^{(m,n)}:=\begin{pmatrix}
F_{2,(k_1,0)}^{(m,n)}&F_{2,(k_1,1)}^{(m,n)}&\cdots&F_{2,(k_1,N)}^{(m,n)}
\end{pmatrix},\\
&\boldsymbol{G}_{2}^{(m,n)}:=\begin{pmatrix}
\boldsymbol{G}_{2,(1)}^{(m,n)}&\boldsymbol{G}_{2,(2)}^{(m,n)}&\cdots&\boldsymbol{G}_{2,(M)}^{(m,n)}
\end{pmatrix}, \quad
\boldsymbol{G}_{2,(k_1)}^{(m,n)}:=\begin{pmatrix}
G_{2,(k_1,0)}^{(m,n)}&G_{2,(k_1,1)}^{(m,n)}&\cdots&G_{2,(k_1,N)}^{(m,n)}
\end{pmatrix},\\
&\boldsymbol{H}_{2}^{(m,n)}:=\begin{pmatrix}
\boldsymbol{H}_{2,(0)}^{(m,n)}&\boldsymbol{H}_{2,(1)}^{(m,n)}&\cdots&\boldsymbol{H}_{2,(M)}^{(m,n)}
\end{pmatrix}, \quad
\boldsymbol{H}_{2,(k_1)}^{(m,n)}:=\begin{pmatrix}
H_{2,(k_1,1)}^{(m,n)}&H_{2,(k_1,2)}^{(m,n)}&\cdots&H_{2,(k_1,N)}^{(m,n)}
\end{pmatrix},
\end{align*}
and
\begin{align*}
F_{2,(k)}^{(m,n)}:=\frac{h}{c^{(m,n)}}2\kappa_0^2\tilde{F}_{2,(k)}^{(m,n)},\quad
G_{2,(k)}^{(m,n)}:=\frac{h}{c^{(m,n)}}\frac{2k_1\pi}{a}\tilde{G}_{2,(k)}^{(m,n)}
,\quad
H_{2,(k)}^{(m,n)}:=\frac{h}{c^{(m,n)}}\frac{-2k_2\pi}{b}\tilde{H}_{2,(k)}^{(m,n)
},
\end{align*}
with
\begin{align*}
&\tilde{F}_{2,(k)}^{(m,n)}:=\int_{\Gamma}\sin\Big{(}\frac{m\pi
x_1}{a}\Big{)}\cos\Big{(}\frac{n\pi x_2}{b}\Big{)}\Big{(}\int_{\Gamma}{\rm
sin}\Big{(}\frac{k_1\pi y_1}{a}\Big{)}\cos\Big{(}\frac{k_2\pi
y_2}{b}\Big{)}g(\boldsymbol{x},\boldsymbol{y}){\rm
d}s_{\boldsymbol{y}}\Big{)}{\rm d}s_{\boldsymbol{x}},\\
&\tilde{G}_{2,(k)}^{(m,n)}:=\int_{\Gamma}\sin\Big{(}\frac{m\pi
x_1}{a}\Big{)}\cos\Big{(}\frac{n\pi x_2}{b}\Big{)}\Big{(}\int_{\Gamma}
\cos\Big{(}\frac{k_1\pi y_1}{a}\Big{)}\cos\Big{(}\frac{k_2\pi
y_2}{b}\Big{)}\partial_{x_1}g(\boldsymbol{x},\boldsymbol{y}){\rm
d}s_{\boldsymbol{y}}\Big{)}{\rm d}s_{\boldsymbol{x}},\\
&\tilde{H}_{2,(k)}^{(m,n)}:=\int_{\Gamma}\sin\Big{(}\frac{m\pi
x_1}{a}\Big{)}\cos\Big{(}\frac{n\pi x_2}{b}\Big{)}\Big{(}\int_{\Gamma}
\cos\Big{(}\frac{k_1\pi y_1}{a}\Big{)}\cos\Big{(}\frac{k_2\pi
y_2}{b}\Big{)}\partial_{x_1}g(\boldsymbol{x},\boldsymbol{y}){\rm
d}s_{\boldsymbol{y}}\Big{)}{\rm d}s_{\boldsymbol{x}}.
\end{align*}
Here
\[
c^{(m,n)}=
\begin{cases}
\frac{ab}{2}, &\mbox{ if } n=0,\\
\frac{ab}{4}, &\mbox{ others. }
\end{cases}
\]

The definition of $\boldsymbol{g}_2$ is given by
\begin{align*}
\boldsymbol{g}_2:=\begin{pmatrix}
\boldsymbol{g}_{2}^{(1)}&\boldsymbol{g}_{2}^{(2)}&
\hdots&\boldsymbol{g}_{2}^{(M)}
\end{pmatrix},\quad
 \boldsymbol{g}_{2}^{(m)}:=\begin{pmatrix}
  g_{1}^{(m,0)}&g_{1}^{(m,1)}&\hdots&g_{1}^{(m,N)}
  \end{pmatrix},
\end{align*}
where
\begin{align*}
g_2^{(m,n)}:=\frac{h}{c^{(m,n)}}2(i\alpha_2p_3+i\beta p_2)\int_{\Gamma}
\sin\Big{(}\frac{m\pi x_1}{a}\Big{)}\cos\Big{(}\frac{n\pi
x_2}{b}\Big{)}e^{i(\alpha_1x_1+\alpha_2x_2)}{\rm d}s_{\boldsymbol{x}}.
\end{align*}

\end{document}